\documentclass{article}
\usepackage{authblk}
\usepackage{amsmath}
\usepackage{amsfonts}
\usepackage{amssymb}
\usepackage{amscd}
\usepackage{mathtools}
\usepackage{gensymb}
\usepackage{graphicx,color}
\usepackage[pdftex]{hyperref}
\usepackage{fullpage}
\usepackage{tikz}
\usepackage{subcaption}
\usepackage{color,soul}
\usepackage{algorithm2e}
\newcommand{\norm}[1]{\left\lVert#1\right\rVert}
\newcommand{\bs}{\boldsymbol}
\newcommand{\ms}{\textsf}

\newcommand{\tb}{\textbf}
\newcommand{\ed}{\text{d}}
\newcommand{\tr}{\text{Tr}}

\title{A novel four-field mixed variational approach to Kirchhoff rods implemented with finite element exterior calculus}
\author[1]{Jamun Kumar N}
\author[2]{Bensingh Dhas}
\author[2]{Arun R Srinivasa\thanks{asrinivasa@tamu.edu}}
\author[2]{J N Reddy\thanks{jnreddy@tamu.edu}}
\author[1]{Debasish Roy\thanks{royd@iisc.ac.in}}
\affil[1]{Centre of Excellence in Advanced Mechanics of Materials, Indian Institute of Science, Bangalore, India}
\affil[2]{J. Mike Walker’66 Department of Mechanical Engineering, Texas A\&M University, College Station, TX 77843-3123, USA}
\date{}
\begin{document}
\maketitle
\begin{abstract}
  A four-field mixed variational principle
  is proposed for large deformation analysis of Kirchhoff rods with the lowest order $C^0$ mixed FE approximations. The
  core idea behind the approach is to introduce a one-parameter family of points (the centerline) and a separate one-parameter family of orthonormal frames (the Cartan moving frame) that are specified independently. The curvature and torsion of the curve are then related to the relative rotation of neighboring frames. The  relationship between the frame and the centerline is then enforced at the solution step using a Lagrange multiplier (which plays the role of section force). 
  Well known frames like the Frenet-Serret are defined only using the centerline,  which demands higher-order smoothness for the centerline approximation. Decoupling the frame from the position vector of the base curve leads to a description of torsion and curvature that is independent of the position information, thus allowing for simpler interpolations. This approach is then cast in the language of exterior calculus. In this framework,  the strain energy may be interpreted as a differential form leading to the natural force-displacement (velocity) pairing. The four-field mixed variational principle we  propose has frame, differential forms, and position vector as input arguments.  While the position vector is interpolated linearly, the frames are interpolated as piecewise geodesics on the rotation group.
  Similarly, consistent interpolation schemes are
  also adopted to obtain finite dimensional approximations for the differential
  forms appearing in the mixed functional. Using these discrete approximations, a discrete mixed variational principle is laid out which is then numerically  extremized. The discrete approximation is then applied to a few benchmark problems. Vis-\'a-vis most available approaches, our numerical studies reveal an impressive performance of the proposed method without numerical instabilities or locking.
\end{abstract}


\section{Introduction}
\subsection{Background}
Rods are often treated as 1-dimensional continuum models where one of the  dimensions is relatively large vis-\'a-vis the rest. Because of the 1-dimensional character, these models bring about computational efficiency and lend themselves to myriad applications, from  DNA supercoiling to computer graphics \cite{spillmann2007corde}. Accordingly, substantive effort has gone into the development of rod
theories and numerical techniques to solve them. Kirchhoff  rods constitute a finite deformation model to study problems with large flexural
deformation. Initially developed by Euler to study special
problems, the model was subsequently refined by Kirchhoff and many others. During
the early days of development, emphasis was placed on linear theories (such as the linearized Euler-Bernoulli beam theory), since
these theories could be used to obtain closed-form solutions for a large
class of problems of engineering interest. Reissner \cite{reissner1972one} extended the Kirchhoff rod model by incorporating shear deformation.
As with Kirchhoff rods, Reissner's shear deformable theory also utilized
stain measures that are frame indifferent and exact up to the deformation
modes assumed in the kinematics. Because of these properties, Kirchhoff and
Reissner theories of rods are generally referred to as geometrically exact
\cite{simo1986three}. Green \textit{et al.} \cite{green1974theory} constructed
a general thermodynamic rod theory based on the idea of a Cosserat curve. By a Cosserat
curve, we mean a curve in  three dimensions with directors or frames attached to
each point on the curve. Since the directors at each point on the curve take
values in $\mathbb{R}^3$, the set of directors at a point may be
identified with the general linear group. Moreover, it is also possible to
reduce Green's rod theory to a Kirchhoff or Reissner rod by incorporating certain constraints.
 
Finite element (FE) approximations for large deformation analysis of rods have taken two
major  directions, one based on Reissner's rod theory and the other using
Kirchhoff rod theory. Even though Kirchhoff rods were used initially, Reissner 
rods that admit shear deformation of the cross section \cite{meier2019geometrically}, 
were favored later on; this choice was dictated by the higher order
smoothness required by the finite element space used to approximate Kirchhoff
rods. In a finite deformation setting, both  Kirchhoff and Reissner rods
can be understood as equations describing a Cosserat or framed curve
\cite{antman2005nonlinear}. However, \emph{the geometric constraint that the frame
should be adapted to the curve, so that one of the vectors in the frame is aligned tangential to the rod's centerline,}\footnote{The Frennet-Serret frame has an even stronger version of adaptability where both the normal and bi-normal vectors are also computed from the curve.} makes the FE approximation of Kirchhoff rods quite challenging. 

Numerical techniques to approximate the solution of a finite deformation Reissner rod was 
introduced by Simo and Vu-Quoc \cite{simo1986three}. One
may recall that Reissner's theory does not incorporate a frame that is
adapted to the rod's base curve.  This FE approximation  extensively used the
rotation group to represent the frame. An important  aspect of their numerical
technique was  the use of the exponential map on $SO(3)$ to update the
configuration variables associated with the frame. This numerical approximation
suffered from locking and issues related to frame indifference. Jelenic
and Crisfield \cite{jelenic1999geometrically}  emphasized the importance of
objectivity and path-independence in the numerical approximation of Reissner
rods. Even though Reissner's rod theory as such is both frame indifferent and path
independent, several numerical techniques used to construct approximate
solutions were shown to be non-objective and path-dependent.  Jelenic and
Crisfield \cite{crisfield1999objectivity} attributed this loss of invariance
to the interpolation scheme used for  approximating the rotation group. They
presented a simple interpolation technique on the rotation group, which could
achieve both path independence and objectivity.  

Approximate solution for Kirchhoff  rods has regained attention in recent
days. These methods are largely based on $C^1$-approximation for deformation. 
In a geometrically linear setting,  Armero and  Valverde \cite{armero2012invariant}
developed a family of FE approximations using a combination of $C^0$ and $C^1$
finite elements. They reported that certain combinations of these finite elements
led to a poor performance of discrete approximations and loss of objectivity. A
similar approach was pursued by Silva \textit{et al.} \cite{e2019simple} in a
large deformation setting, where Hermite interpolation was used to approximate
the position vector and linear approximation was used for the
incremental twist angle.  Meier \textit{et al.} \cite{meier2014objective}
presented an FE approximation for Kirchhoff rods using Hermite interpolation
of displacement fields, which could preserve both path independence and
objectivity. Recently, methods based on isogeometric approximation have gained
prominence as tools to construct numerical approximations for Kirchhoff rods;
these methods offer higher-order continuity across the element boundary. Greco and Cuomo \cite{greco2013b},
presented a purely displacement based isogeometric analysis of Kirchhoff rods in
a linearized setting. This approach was later extended to large deformation
problems by Bauer \textit{et al.}\cite{bauer2016nonlinear}. A major
difficulty with isogemetric analysis for Kirchhoff rods is the requirement of
multiple patches and suitable continuity conditions to enforce the required
smoothness across the patch boundaries.  Moreover, ensuring inter-member
continuity in an assembly of rods is also a challenge since multiple patches have
to be tied together at the member boundaries using additional constraints.
In a restricted setting, Schulz  and B{\"o}l \cite{schulz2019finite} proposed an
FE approximation for  Kirchhoff rods requiring  only $C^0$ continuity.
They incorporated the Kirchhoff constraint  directly into the kinematics,
thereby reducing the independent degrees of freedom to just extensional strain and orientation of the cross section. Their solution procedure essentially had two
steps: the first step predicted the extensional strain and orientation of the
cross section, while the second obtained the deformation by integrating the
differential of position using the previously computed extensional strain and cross-sectional
orientation with Dirichlet boundary condition as
initial values. Algorithmically, this two-step process was carried out using an
additional inner assembly. A drawback of this FE scheme is that
it may not be possible to solve framed structures and problems with multiple
Dirichlet conditions (like a simply supported or continuous beam), since the
integration of displacement (using the  frame and strain) can accommodate only one Dirichlet condition (for each component of deformation). 

Considerable work is also available on modelling Kirchhoff rods from a direct 
discrete standpoint; i.e. these approaches completely sidestep the continuum balance laws and instead directly use Lagrangian mechanics to determine the dynamics of the discrete degrees of freedom. This line of research has its origin in the computer graphics community, where Kirchhoff rods have been used to simulate the
dynamics of slender objects like thread \cite{wang2017real}, human hair \cite{bertails2006super} \cite{daviet2011hybrid}, and so on. Bergou \textit{et al.} \cite{bergou2008discrete}
presented a discrete approximation for Kirchhoff rods, using the idea of
discrete framed curves. In this approach, the discrete configuration of a
Kirchhoff rod is obtained using a simplicial approximation for the base curve. The frame vector tangent to
the discrete curve is defined along the vector connecting the vertices of a
1-simplex. The orientation of the other two vectors is then determined by
the twist angle, an additional degree of freedom on each 1-simplex. A similar effort in this direction is due to Bertails \textit{et al.} \cite{bertails2006super}, who again used Kirchhoff  rods to simulate the dynamics of human hair. In their discrete model, Bertails \textit{et al.} used curvature as the primary variable to describe the configuration of the discrete rod. Using Bishop's frame and the assumption of constant curvature within a simplex, they arrived at a simple expression relating the frame, its curvature and the position vector of the configuration. An important assumption, which permitted the curvature integration is the inextensibilty of the rod. Even though the inextensibility assumption permitted for a simple implementation, it removed the applicability of Bertails \textit{et al.}'s discrete Kirchhoff model to engineering applications like fixed-fixed rods, rods with multiple displacement conditions, and framed structures assembled using multiple rods. Recently, Lestringant \textit{et al.} \cite{lestringant2020discrete} have extended discrete Kirchhoff rods to include inelastic response by introducing inelastic strain as an additional discrete internal variable. 

For Kirchhoff rods, mixed models may be a suitable choice to construct approximations that are both lower order and accurate. Efforts in this direction were taken by Reddy \textit{et al.} \cite{reddy2020nonlinear}, motivated by the finite volume method. They introduced the dual mesh finite domain method for the solution of nonlinear von K\'arm\'an  beams. The key postulate behind the discrete approximation presented by Reddy \textit{et al.} is that the kinetic variables live on the dual cell. This approximation uses finite elements of order at most one to construct the approximate solution. A limitation of Reddy \textit{et al.} is that the beam model cannot be used for general large deformation analysis since the strain measures are not geometrically exact.

\subsection{Motivation to use Cartan's moving frame for rod theory}
In all the works noted above, \emph{the rod's cross-section is used as an a-priori basis for the choice of frame}. Thus, for a Kirchhoff rod, the cross-section remains normal to the centerline and hence requires an adapted frame. On the other hand, for a Reissner's rod or its generalization (e.g. the  Cosserat rod), the cross-section can shear and twist, permitting a frame with no kinematic relationship with the rod's centerline.  In contrast to these, our approach is based on a complete severance of a kinematic relationship between the rod's centerline and the frame --- this decoupling is the basis for an exploitation of Cartan's moving frame. This simple yet significant step allows us to (a) describe the geometrical quantities  independent of the centerline and (b) impose the compatibility between the rod's centerline and the frame \emph{a-posteriori}, leading to a more transparent and elegant description for the geometry of Kirchhoff rods. 

The kinematic decoupling we just mentioned has numerical implications as well. For example, the numerical approximation proposed by Simo and Vu-Quoc uses a frame that is defined by the cross-section and incorporates the strain energy associated with the shear deformation of the cross-section. This shear energy along with interpolation of the frame becomes the source of shear locking in the FE approximation presented by Simo and Vu-Quoc. 
In this work, however, we consider
the Cartan frame independent of both the centerline and the rod cross-section. We then go on to treat the frame adaptation as a constraint enforced through a Lagrange multiplier. Thus, by kinematically decoupling the geometry of the Kirchhoff rod from its position vector, we are able to
\begin{itemize}
  \item construct discrete models for Kirchhoff rods with lower order continuity requirements;
  \item eliminate shear locking without the need for ad-hoc reduced integration approaches and
  \item obtain frame indifferent discretizations.
\end{itemize}

In the present approach, we consider the base curve of the rod as a one-dimensional
(piecewise) smooth manifold embedded within the three-dimensional Euclidean space. The
geometric information about the base curve is encoded within quantities like
curvature and torsion. In a conventional approach,
these geometric quantities are calculated using  higher derivatives of
the curve's position vector. An alternative picture is to keep these geometric
quantities as primal and ask if there are curves with the desired geometric
properties. Fortunately, it is easy to construct such  curves in 1-dimension;
however we do not get one but a family of such curves that are related by
rigid  motion \cite{do2016differential}. This alternative perspective
makes the geometry of the curve as primary quantities and perceives the
placement as a secondary variable, which emerges from the
geometric properties. Cartan's method of moving frames \cite{clelland2017frenet} is an effective tool to
describe a curve using its geometric properties. In this technique, the
geometric information about the curve is described using a collection  of vectors or a frame attached to each point on the curve.  On enforcing the frame adaptivity, the differential of the frame now describes the geometric properties of
the base curve (see Figure \ref{fig:discreteCurve}). Building on this kinematic landscape, we describe the kinetic
quantities of the Kirchhoff rod in terms of the dual frame. Since our description has both position vector and frame as key kinematic variables, additional compatibility conditions are required to ensure that the frame is adapted to the base curve. 

\emph{The scheme outlined above should be contrasted with the use of Frenet-Serret frame or the Bishop frame  (which is a-priori related to the curve geometry) or the usual approaches to Cosserat rods, all of which enforce the compatibility of the frame with the curve a-priori. This a-priori adaptation of the frame is possible only when the centerline is sufficiently smooth. This then creates the challenges of interpolation requiring high order smoothness and/or reduced integration and their attendant issues}. In contrast, the current approach completely ignores any a-priori relationship between the frame and the centerline. Instead, it enforces it as a compatibility constraint a-posteriori. Using this approach, we are able to use simple $C^0$ linear interpolation for both displacement and frame, which comes at the cost of having to interpolate many functions. We enforce compatibility between the frame and the centerline using a variational approach. Towards this, we formulate a mixed variational principle that takes the frame, deformation and certain other differential forms (that appear as Lagrange multipliers) as unknowns. This variational scheme  belongs to the class of variational principle proposed by Dhas and Roy \cite{dhas2020geometric} in the context of three-dimensional elastic solids. The major conceptual difference, however, is the evolution of the frame fields; in other words, the geometry of the configuration evolves, a feature not present is Dhas and Roy \cite{dhas2020geometric}.

Building on the four-field mixed  variational principle, we construct an FE approximation that is only $C^0$. A conventional FE discretization for a Kirchhoff rod approximates the geometry of the base curve from its position vectors. This scheme places a lot of constraints on the approximation space making it difficult to construct. The major idea we put fourth in this paper is to decouple the discrete approximation to the base curve from its geometry.  Since we have the position vector and frame as separate entities, a discrete approximation for the geometry may be effectively performed using the frame, avoiding higher order derivatives. This brings forth two challenges; the first is that the frame fields should satisfy the condition that they are mutually orthonormal -- a nonlinear constraint. The second challenge is to find suitable approximation spaces for the differential forms appearing in the variational principle. We circumvent the first challenge by identifying the frame fields with rotation group $SO(3)$, which incorporates the constraint without the need for a Lagrange multiplier. The problem of finding a discrete approximation to the frame fields is thus solved by using the idea of geodesic interpolation; this idea has its roots in the work of Jeleni\'c and Crisfield \cite{jelenic1999geometrically}. The problem of approximating differential forms is addressed using ideas from finite element exterior calculus (FEEC). Since the problem is 1-dimensional, key features of the FEEC do not appear prominently.
These discrete approximations are then used to construct the discrete mixed functional. Since we have $SO(3)$ as one of the configuration variables, the finite dimensional space on which the solution of the discetre extremization problem is sought becomes a smooth manifold with its own geometric properties. We adopt Newton's method to construct the successive approximations which extremize the discrete mixed functional. Since Newtons's method requires the Hessian of the  functional, we construct it by taking into account the geometry of the discrete configuration space \cite{absil2009optimization}.


The rest of the article is organised as follows. Section \ref{sec:kinematics} introduces the kinematics of a Kirchhoff  rod using the machinery of moving frames. The frame and co-frame fields for reference and deformed configurations of the rod are introduced followed by the  differential of the frame and the connection matrices. The components of the connection matrix are then interpreted  as the curvature and torsion of the base curve. The importance of Cartan's structure equations and its relationship with the underlying geometric assumptions are also highlighted. In Section \ref{sec:HWvarPrin}, various kinematic and kinetic quantities are combined using a  mixed variational principle.  Section \ref{section6:Finite element approximation} discusses a discrete approximation to the four-field mixed variational principle presented in Section \ref{sec:HWvarPrin}. This section also contains the finite dimensional approximations used for different fields appearing in the variational principle. Detailed calculations describing the residual and tangent operators are presented in Appendix \ref{sec:residueTangent}. Using the discrete approximation considered in Section \ref{section6:Finite element approximation}, elaborate numerical experiments are conducted to evaluate the performance of the discrete approximation; these results are presented in Section \ref{sec:NumericalResults}. Finally, Section \ref{sec:conclusion} concludes with a summary. 

\section{Moving frames and kinematics of Kirchhoff rods}
\label{sec:kinematics}
\begin{figure}
  \centering
  \includegraphics{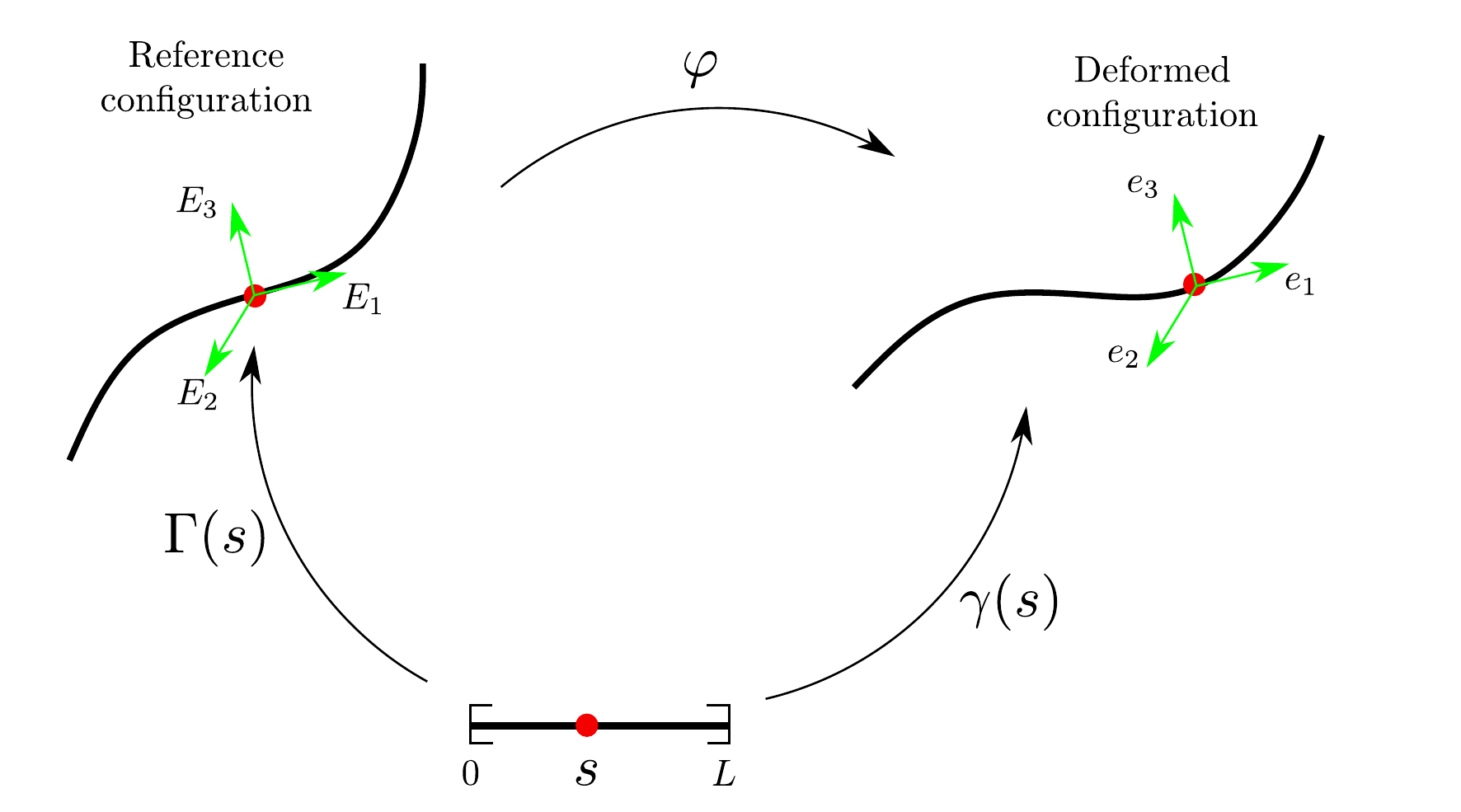}
  \vspace{0.2in}
  \caption{The centerlines of the reference and deformed configuration along with their frame fields}
  \label{fig:referenceAndDefConfigWithFrame}
\end{figure}

In this section, we discuss the kinematics of a Kirchhoff  rod using
Cartan's method of moving frames. The centerline (base curve) of the beam at its reference and deformed configuration are denoted by $\bs{\Gamma}$ and $\bs{\gamma}$ respectively; these are 1-dimensional piece-wise smooth manifolds, with isolated singular
points\footnote{Points where the tangent to the curve is not well defined.}. We also assume that $\bs{\Gamma}$ has an arc-length parametrization so that the centerline of the rod's reference configuration can be
written as $\bs\Gamma:[0,L]\rightarrow\mathbb{R}^3$, where $L$ is the length of the
beam in the reference configuration. Similarly, the centerline of the deformed configuration
can be written as $\bs\gamma:[0,L]\rightarrow\mathbb{R}^3$; however this
map need not be arc-length parametrized (see Fig. \ref{fig:referenceAndDefConfigWithFrame}). The deformation map relates the centerlines of the rod's reference and deformed configurations; we denote it by $\bs{\chi}:\bs{\Gamma}\rightarrow\bs{\gamma}$. This deformation map is more restrictive than the one encountered in 3D continuum mechanics, since the range of $\bs{\chi}$ has co-dimension 2. For calculating the local deformation state of the rod, we use the maps
$\bs{\Gamma}$ and $\bs{\gamma}$, rather than the deformation map, since both $\bs{\Gamma}$ and $\bs{\gamma}$ have the interval $[0,L]$ as their domains. We denote the tangent spaces to the curves $\bs{\Gamma}$  and $\bs{\gamma}$ (at $\tau\in [0,L]$) by $T_{\bs{\Gamma}(\tau)}\bs{\Gamma}$ and $T_{\bs{\gamma}(\tau)}\bs{\gamma}$ respectively, which are vector spaces of dimension one. To describe the position vector of the reference and deformed configurations, we introduce a globally constant frame\footnote{By a frame we mean a collection of vectors orthonormal to each other with respect to the Euclidean inner-product.}, and denote it by $\bar{\mathcal{F}}=\{\bar{\tb{E}}_1,\ldots,\bar{\tb{E}}_3\}$. A natural choice for this frame is the canonical Cartesian basis. In terms of $\bar{\mathcal{F}}$, the differential of the position vector in the reference configuration can be written as,
\begin{equation}
  \ed \tb{X}=\bar{\textbf{E}}_i \frac{d \Gamma^i}{d\tau}d\tau.
  \label{eq:dXFromDeform}
\end{equation}
In the above equation, $\Gamma^i$ denotes the components of the reference configuration's position vector along the direction $\bar{\tb{E}}_i$ and $\frac{d(\cdot)}{d\tau}$ denotes the derivative with respect to the reference configurations arc-length parameter.  
Similarly, the  differential of position in the deformed configuration can be written as,
\begin{equation}
  \ed \bs{x}=\bar{\textbf{E}}_i\frac{d \gamma^i}{d\tau}d\tau,
  \label{eq:dxFromDeform}
\end{equation}
where $\gamma^i$ denotes the components of the deformed configurations position vector. In both \eqref{eq:dXFromDeform} and \eqref{eq:dxFromDeform}, we have retained $d\tau$ on the right hand side. This notation is to highlight the fact that both $\frac{d\Gamma^i}{d\tau} d\tau$ and $\frac{d\gamma^i}{d\tau} d\tau$ are differential forms  on the 1-dimensional manifold $[0,L]$. The map $\ed \tb{X}$ takes a vector from the tangent space of $[0,L]$ and produces a vector from $T_{\bs{\Gamma}(\tau)}\bs{\Gamma}(\tau)$. Similarly, $\ed \bs{x}$ returns a vector from $T_{\bs{\gamma}(\tau)}\bs{\gamma}(\tau)$. An important point here is that, even though both $T_{\bs{\Gamma}(\tau)}\bs{\Gamma}(\tau)$ and $T_{\bs{\gamma}(\tau)}\bs{\gamma}(\tau)$ are both 1-dimensional, we have the derivatives of all the three components appearing in the differential, since we have expanded the basis vectors for $T_{\bs{\Gamma}(\tau)}\bs{\Gamma}(\tau)$ and $T_{\bs{\gamma}(\tau)}\bs{\gamma}(\tau)$  in terms of the Canonical frame $\bar{\mathcal{F}}$.
\subsection{Geometry  via moving frames: The Cartan approach }
In a conventional setting, geometric information like the curvature and torsion of the curves $\bs\Gamma$ and $\bs\gamma$ are obtained by first constructing a frame that is adapted to the curve; the Frennet-Serret frame is commonly used for this purpose. For a curve in 3-dimensions, the Frennet-Serret frame consists of the tangent, normal, and bi-normal vectors. Such a construction demands that the curve be $C^3$ continuous, which is a very strong requirement. Bishop's frame is an alternative to the Frenet-Serret frame; it relaxes the smoothness requirement on the curve by decoupling torsion and curvature. However, incorporating Bishop's frame \cite{bishop1975there} within a structural mechanics application is involved since the torsion of the frame has to be obtained by solving a first order ordinary differential equation along the length of the curve. Similar to Bishop's frame, multiple hybrid frames have been introduced in the literature each adapted for the application at hand \cite{arbind2018higher}. 

In this work, we attach an \emph{independent frame}  at a material point on the rod's base curve. This frame together with its co-frame encode the geometric information such as curvature and torsion that is used for the construction of the energy function of the continuum.  This technique of using a frame and a co-frame to describe the geometry is  called Cartan's method of moving frames. We denote the collection of orthonormal vector fields $\{\tb{E}_1,\ldots,\tb{E}_3\}$ by $\mathcal{F}_{\bs{\Gamma}}$ and call it the frame for the reference configuration. Similarly, the frame of the deformed configuration is denoted by $\mathcal{F}_{\bs{\gamma}}=\{\tb{e}_1,\ldots,\tb{e}_3\}$. Similar to $\tb{E}_i$, the vectors $\tb{e}_i$ are also orthonormal with respect to the Euclidean inner-product. Note that both $\mathcal{F}_{\bs{\Gamma}}$ and $\mathcal{F}_{\bs{\gamma}}$ form a basis for $\mathbb{R}^3$. Since both $\mathcal{F}_{\bs{\Gamma}}$ and $\mathcal{F}_{\bs{\gamma}}$ are orthonormal frames, the relationship between them can be written as,
\begin{equation}
    \tb{e}_i=\bs{\Lambda} \tb{E}_i
    \label{eq:refToDefFrame}
\end{equation}
where $\bs\Lambda$ is an element from the rotation group $SO(3)$. Similarly, the relationship between the globally constant frame $\bar{\mathcal{F}}$ and the frame of the reference configuration can be written as,
\begin{equation}
  \tb{E}_i=\bs{\Lambda}_0\bar{\tb{E}}_i
  \label{eq:fixedToRefFrame}
\end{equation}
Even though the frame $\bar{\mathcal{F}}$ is globally constant, the rotation $\bs{\Lambda}_0$ need not be so since the frame for the reference configuration may depend on the position. Combining \eqref{eq:refToDefFrame} and \eqref{eq:fixedToRefFrame}, we get, 
\begin{equation}
  \tb{e}_i=\bs{\Lambda\Lambda}_0\bar{\tb{E}}_i
\end{equation} 
The co-frames associated with reference and deformed configurations are respectively denoted by $\mathcal{F}^*_{\bs\Gamma}=\{\bs{\alpha}^1,\bs{\alpha}^2,\bs{\alpha}^3\}$ and $\mathcal{F}^*_{\bs\gamma}=\{\bs{\beta}^1,\bs{\beta}^2,\bs{\beta}^3\}$. These co-frames span the cotangent space of $\mathbb{R}^3$ and satisfy the relations $\bs\alpha^i(\tb{E}_j)=\delta^i_j$ and $\bs{\beta}^i(\tb{e}_j)=\delta^i_j$, where $\bs{\alpha}^i(.)$ and $\bs{\beta}^i(.)$ denote the action of co-vectors on tangent vectors. The differential of frame fields in the reference configurations can be written as,
\begin{equation}
  \ed \tb{E}_i=\widehat{\bs{\Omega}} \tb{E}_j;\quad \widehat{\bs{\Omega}}=  \begin{bmatrix}
    0&-\bs{\kappa}_1& \bs{\kappa}_2\\
    \bs{\kappa}_1&0 &-\bs{\kappa}_3 \\
    -\bs{\kappa}_2&\bs{\kappa}_3&0 \\
  \end{bmatrix},
  \label{eq:diffFrameReference}
\end{equation}
where $ \bs{\kappa}_1=\ed \tb{E}_2 \cdot \tb{E}_1; \bs{\kappa}_2=\ed \tb{E}_3 \cdot \tb{E}_1;\bs{\kappa}_3=\ed \tb{E}_2 \cdot \tb{E}_3$. $\widehat{\bs{\Omega}}$ is the skew symmetric connection matrix, whose coefficients are 1-forms; $\bs{\kappa}_i$ are in bold face to emphasize this fact. However, since we are working with a 1-dimensional manifold, these 1-forms can be identified with real valued functions. We have also placed a $\widehat{(\cdot)}$ to highlight the skew-symmetry of the connection matrix. Equations similar to \eqref{eq:diffFrameReference}, are common in rigid body dynamics and mechanics of Cossarat continua. In terms of vector cross-product, \eqref{eq:diffFrameReference} can be written as,
\begin{equation}
    \ed \tb{E}_i=\bs{\Omega}\times \tb{E}_i.  
\end{equation}
In the last equation, $\bs{\Omega}$ is the axial vector associated with the skew-symmetric matrix $\widehat{\bs{\Omega}}$. Using the same argument, the frame differential for the deformed configurations may be written as,
\begin{equation}
  \ed \tb{e}_i=\widehat{\bs{\Phi}} \tb{e}_j;\quad  \widehat{\bs{\Phi}}=  \begin{bmatrix}
    0&-\bs{\tau}_1& \bs{\tau}_2\\
    \bs{\tau}_1&0 &-\bs{\tau}_3 \\
    -\bs{\tau}_2&\bs{\tau}_3&0 \\
  \end{bmatrix},
  \label{eq:diffFrameDeformed}
\end{equation}
where $\bs{\tau}_1=\ed \tb{e}_2 \cdot \tb{e}_1; \bs{\tau}_2=\ed \tb{e}_3 \cdot \tb{e}_1; \bs{\tau}_3=\ed \tb{e}_2 \cdot \tb{e}_3$.   The matrix $\widehat{\bs{\Phi}}$ is called the connection matrix of the deformed configuration. The skew-symmetry in the connection matrices is due to the frame orthogonality. The differentials given in \eqref{eq:diffFrameReference} and \eqref{eq:diffFrameDeformed} describe local change to the frames at a material point along the centerline of the rod in the reference and deformed configurations. The connection matrices $\widehat{\bs{\Omega}}$ and $\widehat{\bs{\Phi}}$ encode infinitesimal parallel  transport of vectors on the manifolds $\bs{\Gamma}$ and $\bs{\gamma}$, respectively. However, in keeping with the convention in structural mechanics, we call  $\bs{\kappa}_i$ and $\bs{\tau}_i$, $i=1,2$, the curvatures of reference and deformed configurations; $\bs{\kappa}_3$ and $\bs{\tau}_3$ are called the torsions of the respective configurations. These 1-forms describe the rotation of the frame as one moves along the curves $\bs{\Gamma}$ and $\bs{\gamma}$. 

Until now, we have described the tangent to the curve and properties of the moving frame separately. We now introduce the compatibility between them as a \emph{constraint} to be imposed on the frame. This constraint allows us to conflate the geometry of the frame with that of the underlying curve. 

\subsection{The adaptability constraint}
In the last subsection, we introduced the frames $\mathcal{F}_{\bs{\Gamma}}$ and $\mathcal{F}_{\bs{\gamma}}$ independent of the base curves $\bs\Gamma$ and $\bs\gamma$. 
The adaptability condition constrains the frame fields $\tb{E}_1$ and $\tb{e}_1$ to be vectors tangential to the base curves $\bs\Gamma$ and $\bs\gamma$.  In \eqref{eq:dXFromDeform} and \eqref{eq:dxFromDeform}, we wrote the differential of position in terms of the fixed frame $\bar{\mathcal{F}}$. Adaptable frame also means that the tangent space of the curve is spanned by a subset of the frame fields; in the present case, it is only the first frame field that spans the tangent space of the base curves.  A convenient and economic approach to  describe the differentials of position vector is to describe them in terms of the basis of  $T_{\bs{\Gamma}(\tau)}\bs{\Gamma}(\tau)$ and $T_{\bs{\gamma}(\tau)}\bs{\gamma}(\tau)$ respectively. The adaptivity of the frame to the centerline, i.e. aligning the tangent to the centerline along the frame field $\tb{E}_1$, implies that,
\begin{align}
    \label{eq:dXdiff}
  \ed\tb{X}&=\tb{E}_1\otimes \bs{\xi},\\
  \label{eq:dxdiff}
  \ed\bs{x}&=\tb{e}_1\otimes \bs{\eta},
\end{align}
where the quantities $\bs{\xi}$ and $\bs{\eta}$ are differential forms on the interval $[0,L]$. The relations given in \eqref{eq:dXFromDeform}, \eqref{eq:dxFromDeform}, and \eqref{eq:dXdiff}, \eqref{eq:dxdiff} are exactly the same; however they are written in two different bases. The description in \eqref{eq:dXdiff} and \eqref{eq:dxdiff} uses one vector and a differential form, while the relations given in \eqref{eq:dXFromDeform} and \eqref{eq:dxFromDeform} utilize three vectors and three differential forms. Since we now have a frame adapted to the rod's centerline, the co-tangent spaces of the reference and deformed configurations of the rod is spanned by the co-vectors $\bs{\alpha}^1$ and $\bs{\beta}^1$; we denote these one-dimensional spaces as $T^*_{\bs{\Gamma}(\tau)}$ and $T^*_{\bs{\gamma}(\tau)}$, respectively. Using the frame and co-frame fields, the differential of the deformation can be written as,
\begin{equation}
    \ed \bs{\chi}=\tb{e}_1\otimes\bs{\theta}.
\end{equation}
In the above equation, $\bs{\theta}$ is an element from $T^*_{\bs{\Gamma}(\tau)}\bs{\Gamma}$, it is obtained as the pull-back of $\bs \eta$ under the deformation map. Since $T^*_{\bs{\Gamma}(\tau)}\bs\Gamma$ is just one dimensional, $\bs{\theta}$ will be $c\bs{\alpha}^1$, where $c$ is a scalar. However, this scalar depends on the position vector of the rod. Moreover, this scalar $c$ is what encodes the extension experienced by the deformed configuration. Rigorously, $\bs{\theta}$ is the pull-back of $\bs{\beta}^1$ under the deformation map.  Using the definition of a differential form, it is easy to see that range of $\ed \bs{\chi}$ is $T_{{\bs{\gamma}(\tau)}}\bs{\gamma}$. 

In the context of structural mechanics, the cross-section of the rod provides a natural orientation for the frame fields $\tb{E}_i$ and $\tb{e}_i$, $i=2,3$. These frame fields can be chosen to align along the principle directions\footnote{Or some other orthonormal directions, in which the cross-section moment of inertia is known.}  of the cross-sectional moments of inertia. Choosing such an orientation for the other two frames may not be possible when one uses  conventional frames like Frenet-Serret or Bishop, since the other two frame fields are also determined by the rod's base curve. Moreover, for a Kirchhoff rod, the orientation of the cross-section plays an important role in describing the deformation state of the rod at any material point. The Kirchhoff constraint requires the cross-section of the rod to be aligned normal to the tangent vector $\tb{E}_1$ in reference configuration and $\tb{e}_1$ in the deformed configuration.  To describe the orientation of the rod's cross-section, we introduce the orthonormal vector $\tb{E}_2$ and $\tb{E}_3$ for the reference configuration and $\tb{e}_2$ and $\tb{e}_3$ for the deformed configuration. The reference configuration's cross-section may now be described as a subset of the plane spanned by the vectors $\tb{E}_2$ and $\tb{E}_3$. Similarly, the deformed configuration's cross-section is contained in the plane spanned by the vectors $\tb{e}_2$ and $\tb{e}_3$. Our present description of frames does not require the centerline of the beam to be a $C^3$ curve. For our purposes, any frame whose first frame field is aligned along the centerline's tangent and has the same orientation along the curve will be a valid frame. This freedom is possible since geometric information like the curvature and torsion is encoded not in the frame but in its differentials. An important challenge in using moving frames to describe a configuration (manifold) is the enforcement of compatibility between the frame and its configuration. Generally, for the differential of position and frame to be compatible, one requires Cartan's structure equations to be satisfied \cite{clelland2017frenet}. However, in 1-dimension, these equations are trivially satisfied since the topology is always trivial\footnote{ Integrability in 1-dimension just boils down to ordinary integration.}. Moreover,  since the dimension of the co-tangent space is just one, the components of the connection matrices $\bs{\kappa}_i$ and $\bs{\tau}_i$ in \eqref{eq:diffFrameReference} and \ref{eq:diffFrameDeformed} may be identified with real numbers; this property will be extensively used is in the following sections.

\section{A four-field mixed variational principle}
\label{sec:HWvarPrin}
In Section \ref{sec:kinematics}, we discussed the kinematics of Kirchhoff rods using the method of moving frames. We now present a mixed variational principle that takes deformation, frame fields, section force and deformation 1-form as independent variables. This variational principle has the same spirit as the mixed variational principle for 3D elastic solids discussed in \cite{dhas2020geometric}; for more details about mixed methods in mechanics, we refer to Oden and Reddy \cite{oden2012variational}.  The mixed energy functional for a Kirchhoff  rod is postulated as,
\begin{eqnarray}
    I=\int_{[0,L]}[W_b+W_e-\bs{n} (\tb{e}_1\otimes \bs{\eta}  - \ed \bs{x})-\bs{q}\cdot \bs{\gamma}]d\tau-(\bs{\gamma} \cdot \bs{n})|_{\partial [0,L]}-\langle\widehat{\bs{m}},\widehat{\bs{\Omega}}\rangle|_{\partial [0,L]},
    \label{eq:HWvariationalPrinciple}
\end{eqnarray}
where $W_b$ and $W_e$ denotes the energy due to bending and extension, $\bs{q}$ is the externally applied  force density and the last two terms can be interpreted as  the forces and moments applied at the boundary of the rod.
In the above functional, the  Lagrange multiplier which enforces the adaptability constraint between frame and position vector can be interpreted as the section force $\bs{n}$.   In terms of the co-frames, the force  acting on the deformed cross-section can be written as,
\begin{equation}
  \bs{n}=n_i\bs{\beta}^i
  \label{eq:forceOneForm}
\end{equation}
where $n_i$ are real numbers. Conventionally, forces are identified with vectors. Here, we identify them with differential forms to be consistent with the force-velocity pairing.
\begin{figure}
  \centering
  \includegraphics[scale=0.5]{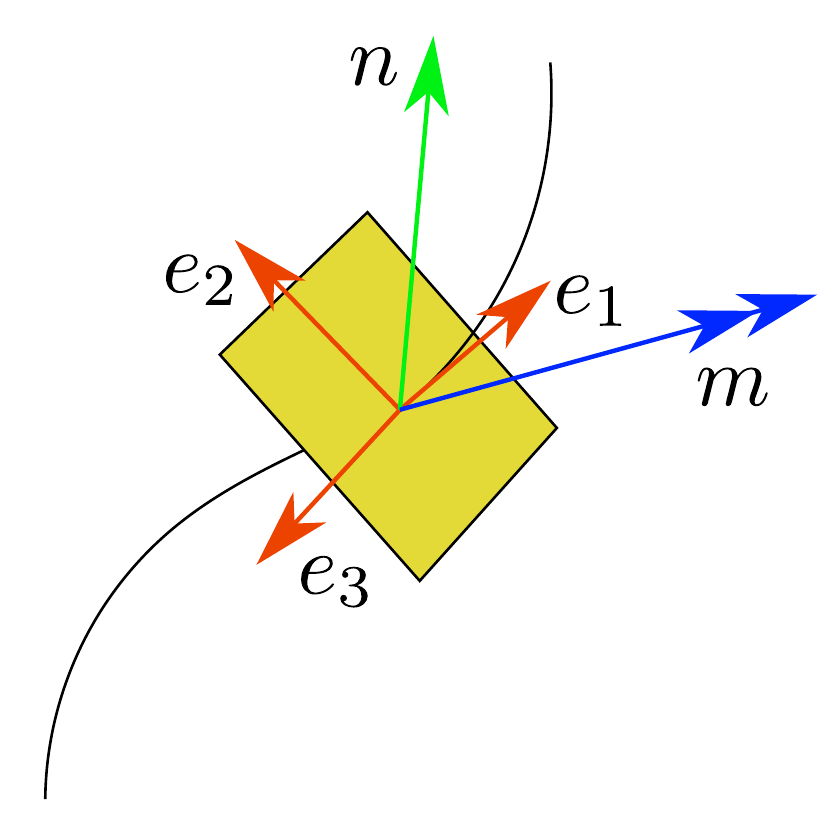}
  \caption{The deformed configurations centerline with a cross-section  showing the resultant force and moment.}
  \label{fig:crossSectionResultants}
\end{figure}
This four-field functional is mixed only in deformation as \emph{moments are not considered independent in the functional}. As already stated, the space $T^*_{\bs{\Gamma}(\tau)}\bs{\Gamma}$ is only 1-dimensional, which permits us to identify these spaces with real numbers. The functional may thus be written as,
\begin{eqnarray}
  I=\int_{[0,L]}\left[W_b+W_e-\left(\eta \bs{n}(\tb{e}_1)-\frac{d \gamma^i}{ds}\bs{n}(\bar{\tb{E}}_i)\right)-\bs{q}\cdot \bs{\gamma}\right]d\tau-(\bs{\gamma}\cdot \bs{n})|_{\partial [0,L]}-\langle\widehat{\bs{m}},\widehat{\bs{\Omega}}\rangle|_{\partial [0,L]}.
\end{eqnarray}
Now using the definition of $\bs{n}$ in terms of the co-frame in the deformed configuration, the above equation becomes,
\begin{equation}
  I=\int_{[0,L]}\left[W_b+W_e-\left( \eta \bs{n}(\tb{e}_1)-\frac{d \gamma^i}{ds}\bar{\bs{\Lambda}}^j_i n_j\right)-\bs{q}\cdot \bs{\gamma}  \right]d\tau -(\bs{\gamma}\cdot \bs{n})|_{\partial [0,L]}-\langle\widehat{\bs{m}},\widehat{\bs{\Omega}}\rangle|_{\partial [0,L]},
  \label{eq:HWvartiationalPrinciple}
\end{equation}
where $\bar{\bs{\Lambda}}^i_j=\bs{\beta}^i(\bar{\tb{E}}_j)$. The third and fourth terms in the mixed functional enforces the constraint that the deformation gradient computed from the deformation map should be identical to the deformation gradient constructed by pulling the differential of the deformed configuration back to the reference configuration. In other words, we want the differential given in \eqref{eq:dxdiff} to be integrable. The section force acts as a Lagrange multiplier enforcing  this integrability  constraint. This integrability condition on the deformation gradient also contains the Kirchhoff constraint. 

\subsection{Bending, twisting and extension energy}
We assume that the strain energy density of the rod is  quadratic in curvatures; the form of the
this  energy is given as,
\begin{equation}
  W_b=\frac{1}{2}\langle \widehat{\bs{\Omega}}, \mathbb{C}:\widehat{\bs{\Omega}} \rangle_{T_I SO(3)}.
  \label{eq:bendingEnergy}
\end{equation}
Here, $\mathbb{C}$ is a forth order constitutive tensor and $:$  denotes contraction between a fourth order tensor and a second order tensor. The constitutive tensor is positive definite on the set of skew symmetric tensors. If the principle directions of the cross section moment of inertia are aligned along the $\bar{\tb{E}}_i$ directions, the constitutive relation between the moments becomes decoupled. For such a constitutive relationship, the action of $\mathbb{C}$ on the basis elements of $T_I SO(3)$ is given by,
\begin{equation}
  \mathbb{C}:\widehat{\bs{w}}_1=c_1\widehat{\bs{w}}_1;\quad\mathbb{C}:\widehat{\bs{w}}_2=c_2\widehat{\bs{w}}_2;\quad\mathbb{C}:\widehat{\bs{w}}_3=c_3\widehat{\bs{w}}_3,
\end{equation}
where $\bs{w}^i$  are the basis vectors for $so(3)$; more details about $so(3)$ can be found in Appendix \ref{section2:Rotation group}.  The constants $c_1$, $c_2$ and $c_3$ can be identified with the bending and torsional stiffnesses of the rod, namely $EI_1$, $EI_2$ and $GJ$.  From the above definition, the following expression can be written for $\mathbb{C}$,
\begin{equation}
  \mathbb{C}=c_1(\widehat{\bs{w}}_1\otimes \widehat{\bs{w}}_1)+c_2(\widehat{\bs{w}}_2\otimes \widehat{\bs{w}}_2)+c_3(\widehat{\bs{w}}_3\otimes \widehat{\bs{w}}_3).
\end{equation}
With this, the derivative of the energy density with respect to curvature can be identified as the moment tensor that may now be written as,
\begin{equation}
  \widehat{\bs{m}}=\frac{\partial W_b}{\partial \widehat{\bs{\Omega}}},
\end{equation}
which is skew-symmetric and identifiable as an element from $T^*_ISO(3)$. Thus the section forces (which appear as constraints and the moments) can be visualized as in  Figure \ref{fig:crossSectionResultants} . 

The constitutive tensor can be thought of as a map sending $T_I SO(3)$ to $T_I^* SO(3)$. 
From the above definition, it is clear that the moment tensor is also skew
symmetric. In terms of the axial vectors of $\widehat{\bs{\Omega}}$ and $\widehat{\bs{m}}$, the constitutive relation can be written as,
\begin{equation}
  \bs{m}=\tb{D}\bs{\Omega}
\end{equation}
The expression for bending energy given in \eqref{eq:bendingEnergy} can be pulled under the hat-mat to get a simpler expression,
\begin{equation}
    W_b=\frac{1}{2}\langle \bs{\Omega}, \tb{D}\bs{\Omega} \rangle_{\mathbb{R}^3}
    \label{eq:consMomentR3}
\end{equation}
In the above expression, $\tb{D}$ is a diagonal matrix with constant $c_1$, $c_2$ and $c_3$ along the diagonal.
Similar to the bending energy, the energy due to extension is also written as a quadratic function. In terms of the deformation 1-form, the  strain energy density due to rod extension can be written as, 
\begin{equation}
  W_e=\frac{AE}{2}(\eta -1)^2,
\end{equation}
Here, $AE$ denotes the extensional stiffness of the rod.

\section{Finite element approximation}
\label{section6:Finite element approximation}
Apart from the clarity brought about by separating the curvature and extensional strains from the centerline properties, a more pronounced advantage of this approach can be seen when we consider finite element approximations.

\begin{figure}
  \centering
 \includegraphics[scale=0.59]{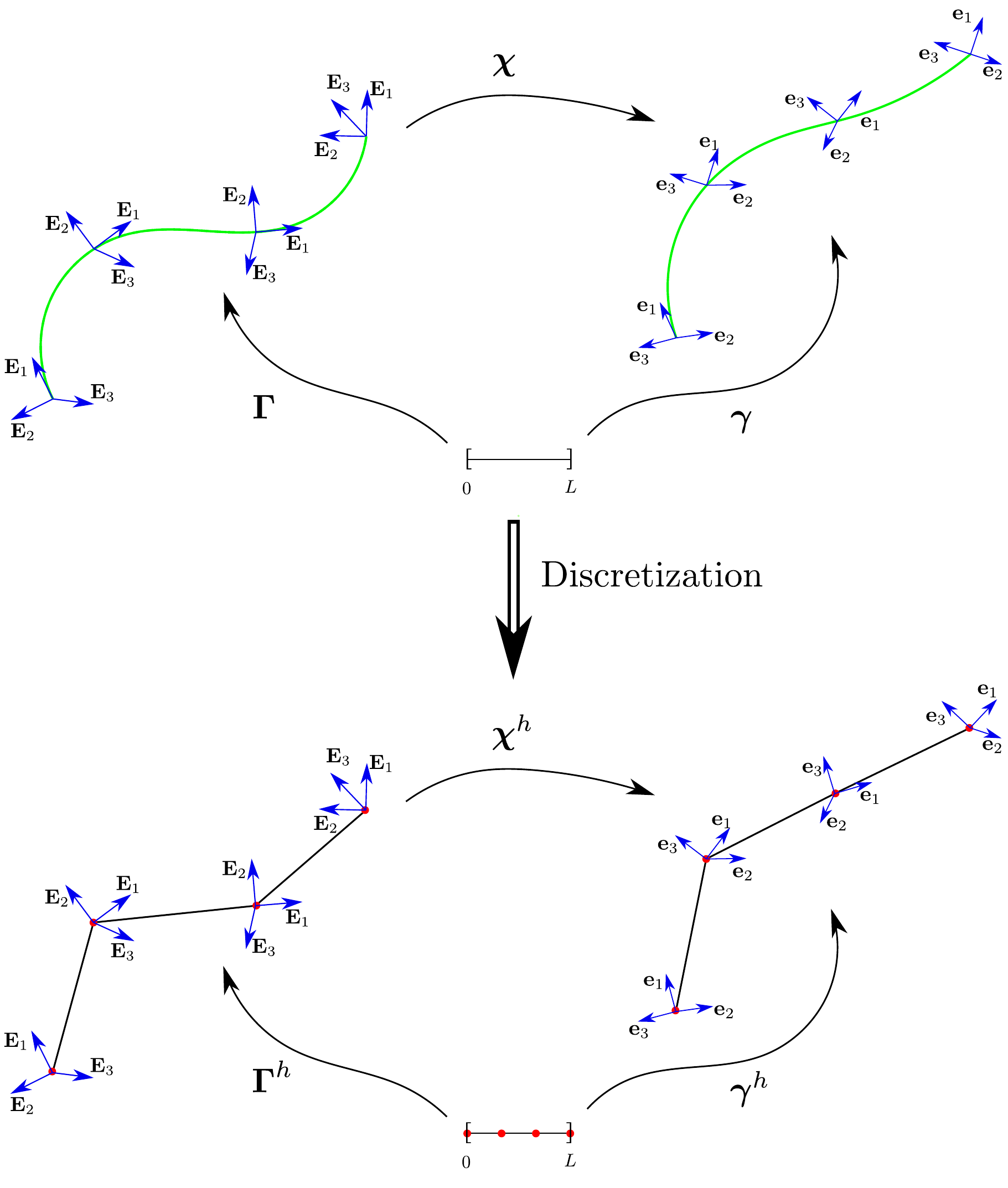}
 \caption{The reference and deformed configurations of the rod along with its discrete approximation. We have not aligned the frame field $\tb{e}_1$ along the deformed configuration's tangent vector, this is to highlight the fact that the frame is chosen independently of the curve.  Note that for the discrete case, for both the reference and deformed configurations, the frame  is not tangential to the edge connecting the vertices.}
 \label{fig:discreteCurve}
\end{figure}

In this section, we present an FE approximation for the four-field mixed variational principle for Kirchhoff rods discussed in Section \ref{sec:HWvarPrin}. The  variational principle has deformation, section forces, deformation 1-forms and frame fields as independent arguments. We seek to independently approximate these quantities using the lowest order finite elements. In Section \ref{sec:kinematics}, the frame fields were identified with the rotation group; this identification will be used in computation as well.  The variational principle has differential operators of degree at most one associated with the frame fields and deformation. The traction and deformation 1-forms do not have any derivatives associated with them. In the present approximation, we use the Lagrange FE space of degree one to approximate the displacement components. The discrete approximation for the $i$th displacement  component can be written as,
\begin{equation}
  u_h^i=\sum_{j=1}^{2}N^ju^i_j;\quad i=1,2,
  \label{eq:discreteApproxDisp}
\end{equation}
where the subscript $h$ indicates that the field is only an approximation, $N^j$ denotes $j$th Lagrange shape functions and $u^i_j$ denotes the $j$th degrees of freedom associated with the $i$th displacement component. The discrete approximation for the differential of deformed position can be written as,
\begin{equation}
  \ed \bs{x}^h=\ed \tb{X}^h+\ed \bs{u}_h,
\end{equation}
where $\ed \bs{x}$ and $\ed \tb{X}$ denote the discrete approximations for the differential of position in reference and deformed configurations. In terms of the frame $\bar{\mathcal{F}}$, the discrete approximation for displacement can be written as,
\begin{equation}
  \bs{u}_h=\sum_{i=1}^3 u_h^i\bar{\tb{E}}_i.
\end{equation}
In writing the above equation, we have used the fact that $\bar{\tb{E}}_i$ are known a-priori and only discrete approximation for the displacement components need to be found. Introducing $\tb{N}:=[N_1,\ldots,N_n]$, $n$ being the number of nodes per element, the discrete approximation for each displacement component can be written as $u^i_h=\tb{N}\tb{u}^i$\footnote{$\tb{N}$ is a row vector and $\tb{u}^i$ is a column vector.}, where $\tb{u}^i$ is the vector obtained by collecting the degrees of freedom associated with the $i$th component of displacement. The differential of displacement can be written as,
\begin{equation}
  \ed \bs{u}_h=\sum_{i=1}^3 \bar{\tb{E}}_i\frac{d u_h^i}{d\tau}d\tau.
\end{equation}
Using the previously introduced notation, the $i$th differential of displacement can be written as $\ed u^i=\ed \tb{N}\tb{u}^i$ with $\ed \tb{N}=[\ed N_1,\ldots,\ed N_n]$.
Similarly, the reference configuration's differential can be written as,
\begin{equation}
  \ed \tb{X}^h=\sum_{i=1}^3\bar{\tb{E}}_i \frac{d \Gamma_h^i}{d\tau}d\tau.
\end{equation}
where $\Gamma_h^i$ denotes the $i$th component of position, described using the frame $\bar{\mathcal{F}}$.
Since the reference configuration is known a-priori, establishing a 1-dimensional finite element mesh on the reference configuration would provide the components $\Gamma^i_h$. Figure \ref{fig:discreteCurve} describes the discrete approximation for centerline curves of the reference and deformed configurations along with the discrete frame fields. 

The discrete approximation for the rotation field is a more involved since any valid  discrete approximation should produce values in the rotation group. Like the displacement field, we seek a $C^0$ approximation for the frame fields as well.  The conditions $\bs{\Lambda}^t\bs{\Lambda}=\tb{I}$ and $\det(\bs{\Lambda})=1$ prevent us for using the classical idea of interpolation for the rotation group directly.
We now present a discrete approximation for the  rotation field using the exponential map on $SO(3)$. Given two rotations $\Lambda_1$ and
$\Lambda_2$, describing the frame fields at each node of a finite element mesh, the
objective  is to find an approximation for the frame fields. In a conventional finite element technique, fields can be directly
interpolated since the interpolated quantity takes value in a vector space. We circumvent the difficulties, as noted above, by using a geodesic interpolation scheme on $SO(3)$. In structural mechanics, these approximation schemes have their origins in the work of Simo and Vu-Quoc \cite{simo1986three}. The
incremental rotation, $\Delta \Lambda$ relating the rotations  $\Lambda_1$ and $\Lambda_2$ is given as,
\begin{equation}
  \bs{\Lambda}_2=\Delta \bs{\Lambda}\bs{\Lambda}_1.
  \label{eq:incerementWithinElement}
\end{equation}
Note that $\Delta \Lambda$ is also a
finite rotation. The above equation can be further written as,
\begin{equation}
  \bs{\Lambda}_2=\exp(\widehat{\bs{\psi}})\bs{\Lambda}_1.
\end{equation}
where $\bs{\psi}$ is a constant vector which needs to be determined.
The relationship between $\widehat{\bs{\psi}}$ and $\Delta\bs{\Lambda}$ is given by the
$\log$ map on $SO(3)$, formally written as,
\begin{equation}
  \widehat{\bs{\psi}}=\log(\Delta \bs{\Lambda}).
\end{equation} 
The incremental rotation can be computed from \eqref{eq:incerementWithinElement} as, $\Delta \bs{\Lambda}=\bs{\Lambda}_2\bs{\Lambda}_1^t$. Using the vector $\bs{\psi}$, the rotation field within the finite element can be written as,
\begin{equation}
  \bs{\Lambda}^h=\exp\left(\frac{s}{L^h}\widehat{\bs{\psi}}\right)\bs{\Lambda}_1;\quad s\in [0,L]
  \label{eq:expInterpolation}
\end{equation}
where $s$ is a parameter along the length of the finite element.
We formally denote the above relationship as a map $\bs{\Lambda}^h:[0,L^h]\rightarrow
SO(3)$, which depends on $\bs{\psi}$ and $\bs{\Lambda}_1$.
In the above equation, a superscript $h$ is added to highlight that
\eqref{eq:expInterpolation} is only an approximation. The relationship given in
\eqref{eq:expInterpolation} is the solution to the differential equation,
\begin{equation}
  \frac{d \bs{\Lambda}^h}{ds}=\frac{\widehat{\bs{\psi}}}{L}\bs{\Lambda}^h,
\end{equation} 
with initial condition $\bs{\Lambda}^h|_{s=0}=\bs{\Lambda}_1$. If one interprets the columns of $\bs{\Lambda}^h$ as the frame fields, the above FE approximation for the rotation field indicates that the connection 1-forms (curvature and torsion in structural mechanics) are constant within the finite element.

\begin{figure}
  \begin{subfigure}{.5\textwidth}
    \centering
    \includegraphics[scale=0.4]{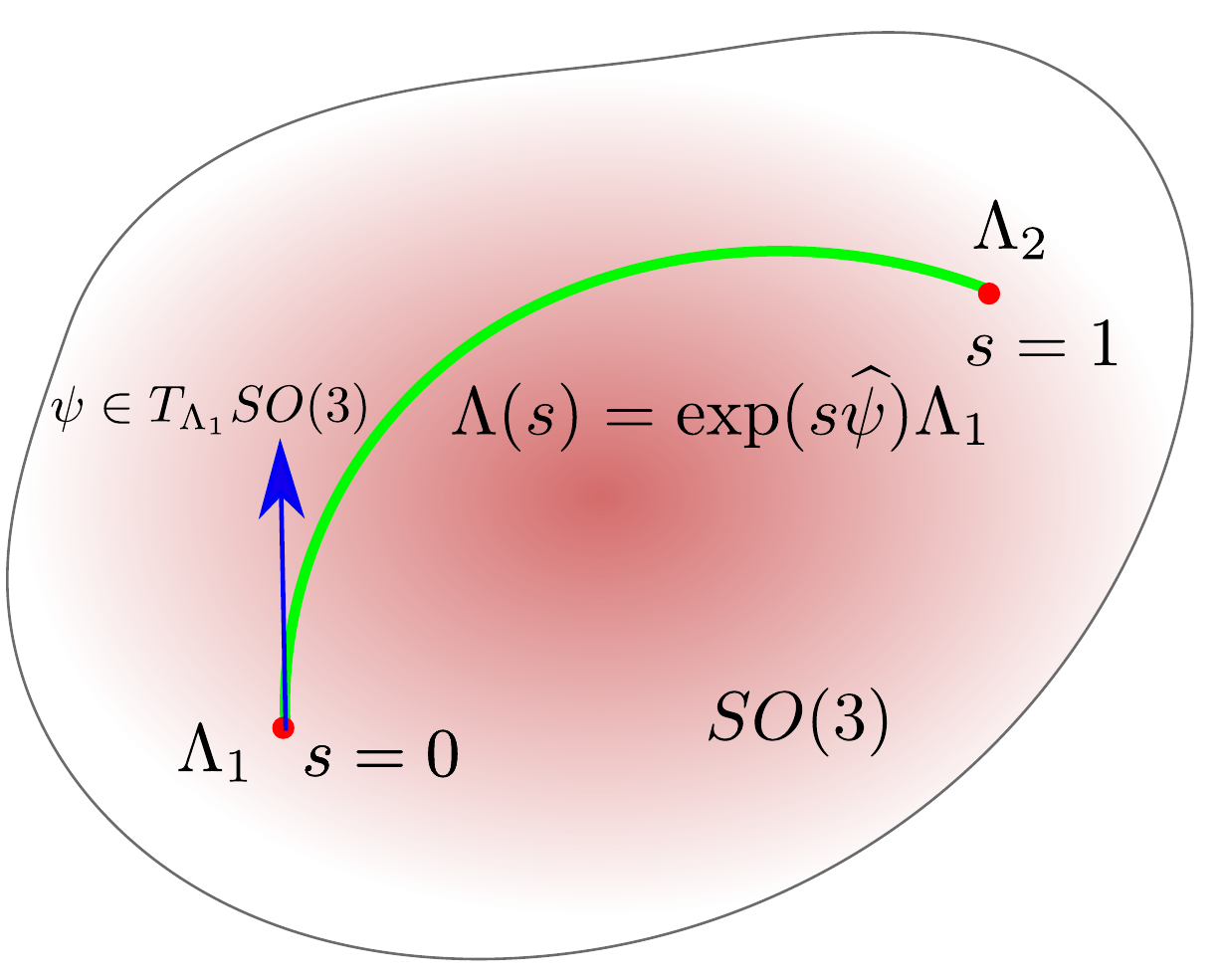}
    \caption{}
  \end{subfigure}%
  \begin{subfigure}{.5\textwidth}
    \centering
    \includegraphics[scale=0.4]{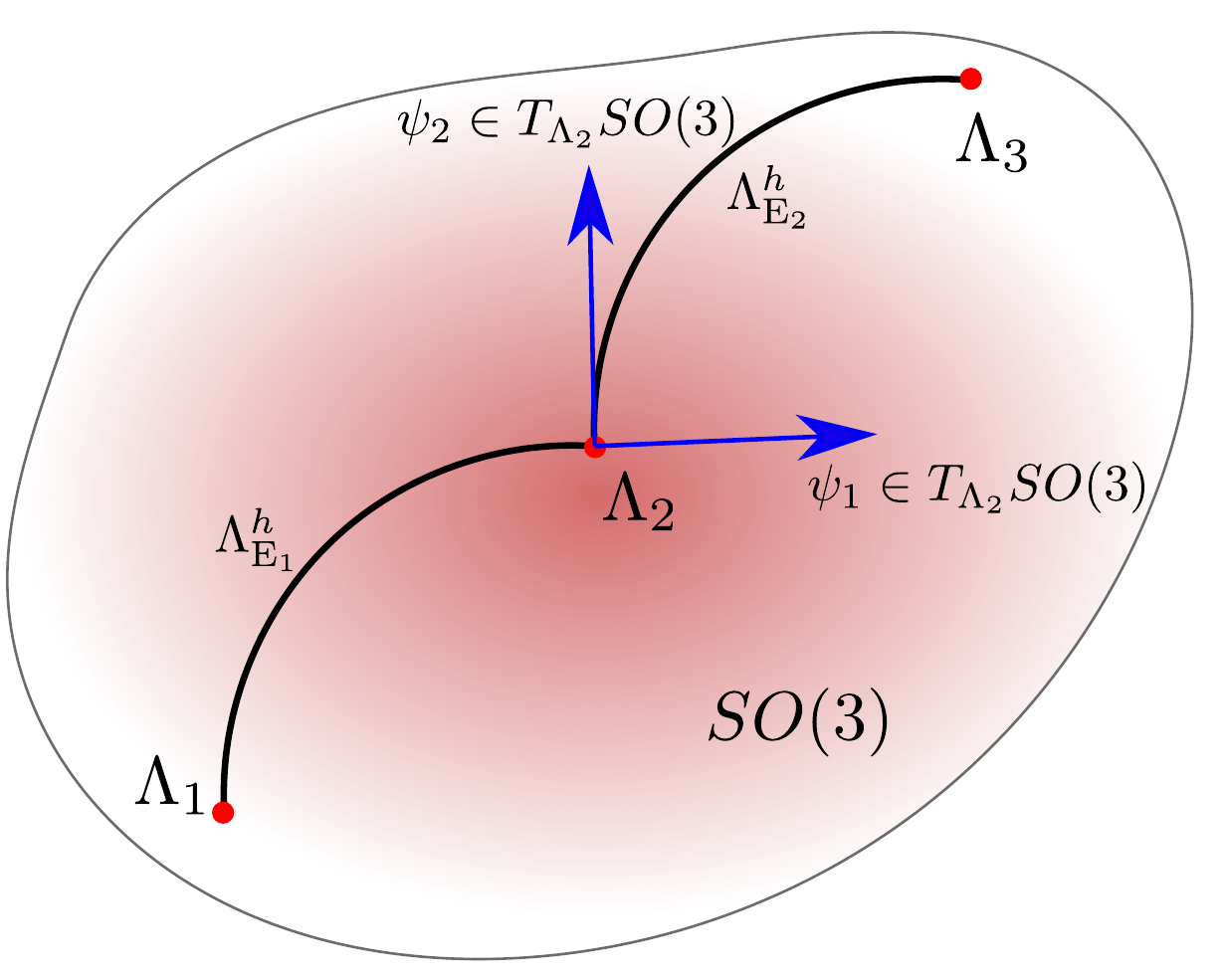}
    \caption{}
  \end{subfigure}
  \caption{(a) Discrete approximation for the frame fields using the rotation group. (b) Blue arrows denote the tangent to the curves $\Lambda_{\text{E}_1}^h$ and  $\Lambda_{\text{E}_2}^h$ at the point $\Lambda_2$; note that the arrows do not overlap, making the approximation only $C^0$.}
  \label{fig:approximationRotation}
  \end{figure}

We now study the continuity properties of the proposed FE approximation for rotation; for this purpose we consider two elements E$_1$ and E$_2$ intersecting at a common node. The approximate rotation fields within E$_1$ and E$_2$ are denoted by $\bs{\Lambda}^h_{\text{E}_1}$ and $\bs{\Lambda}^h_{\text{E}_2}$ respectively; $\bs{\Lambda}^h_{\text{E}_1},\bs{\Lambda}^h_{\text{E}_2}:[0,L^h]\rightarrow SO(3)$. We also have
$\bs{\Lambda}^h_{\text{E}_1}|_{s=0}=\bs{\Lambda}_1$, $\bs{\Lambda}^h_{\text{E}_1}|_{s=L^h}=\bs{\Lambda}^h_{\text{E}_2}|_{s=0}=\bs{\Lambda}_2$ and
$\bs{\Lambda}^h_{\text{E}_2}|_{s=L^h}=\bs{\Lambda}_3$, Figure \ref{fig:approximationRotation}(b) illustrates
these details. We define $\bs{\eta}_1=\frac{d \bs{\Lambda}^h_{\text{E}_1}}{ds}|_{s=L^h}$ and
$\bs{\eta}_2=\frac{d \bs{\Lambda}^h_{\text{E}_2}}{ds}|_{s=0}$. It can be clearly seen that $\bs{\eta}_1,
\bs{\eta}_2 \in T_{\Lambda_2}SO(3)$, however $\bs{\eta}_1$ need not be equal to $\bs{\eta}_2$ since the $\bs{\eta}_1$ and $\bs{\eta}_2$ are tangent to two different curves in $SO(3)$.

Having introduced the approximation for the kinematic fields, we now construct discrete approximations for the force and deformation 1-forms. The discrete approximation for the force 1-form can be written as,
\begin{equation}
  \bs{n}^h=\sum_{i=1}^3n_i^h(\bar{\bs{\Lambda}}^h\bar{\tb{E}}^i),
\end{equation}
where $n_i^h$ are the discrete approximations for the components of the 1-form. Note that in the above equation, we have written the deformed configurational co-frame as a function of the fixed co-frame $\mathcal{F}^*$. The approximate components $n_i^h$ can thus be written as,  
\begin{equation}
  n^h_i=\phi\bar{n}_i
\end{equation}
where, $\phi$ is a function whose value is $1$ within the finite element and zero outside and $\bar{n}_i$ denotes the DOF associated with the $i$th component of force 1-from.  An important point here is that, even though we use constant shape functions to approximate the components of force 1-forms, the force 1-form need not be constant within the finite element since the co-frame fields vary across the length of the element. Our interpolation strategy for the force 1-form is not arbitrary; this choice is motivated by the idea of a simplicial complex and their co-chains. If we consider a simplicial complex as a discrete approximation for a smooth manifold, then the co-chain complex provides an approximation for the differential forms defined on the manifold \cite{hirani2003discrete}. From such a perspective, the discrete approximation for a differential form of degree one can be represented as real numbers defined on the edges. Indeed, the present approximation for force 1-forms is motivated by these ideas. Alternatively, ideas proposed by Arnold and co-workers \cite{arnold2006finite} can be used to construct finite dimensional approximations for differential forms on a simplex; however such elaboration may not be required for a rod which is just a 1-dimensional manifold. 

The final quantity that needs to be discretized is the deformation 1-form. The finite element approximation for the deformation 1-forms can be written as,
\begin{equation}
  \bs{\eta}^h=\phi\bar{\eta}
\end{equation}
where $\bar{\eta}$ denotes the DOF associated with the deformation 1-form. Similar to the force 1-form, we have used the constant FE space to approximate the deformation 1-form. Even though both $\bs{n}$  and $\bs{\eta}$ are differential forms of degree one, the number of DOF's required to describe these quantities are different; this is because $\bs{n}$ is a differential form on $\mathbb{R}^3$, while $\bs{\eta}$ is one on  $\mathbb{R}$. Table \ref{tab:fieldsAndApproximation} gives a summary of all the FE spaces used to discretize different fields and the dimensions of their FE spaces. To summarize, the present FE approximations have DOF's both at nodes and edges. Counting all, a single element would have a total of 16 DOF's.

\begin{table}[h]
  \centering
  \caption{Finite element spaces used to approximate different fields }
  \begin{tabular}{llr}
    \hline
   Field & FE space & \# DoF per element\\
   \hline
   Displacement & linear Lagrange& $2\times 3=6$ \\
   Rotation& linear Lagrange & $2\times 3=6$ \\
   Traction 1-form& constant & $1 \times 3=3$ \\
   Deformation 1-form & constant & 1\\
   \hline
  \end{tabular}
  \label{tab:fieldsAndApproximation}
\end{table}

\subsection{Discrete mixed variational principle}
Having introduced the discrete approximation for frame fields, deformation 1-form, force 1-form and deformation, we now proceed to discuss the discrete approximation for the mixed functional presented in \eqref{eq:HWvartiationalPrinciple}. The discrete functional is obtained by replacing the fields in \eqref{eq:HWvartiationalPrinciple} by their discrete approximations developed in the previous sub-section. The discrete mixed functional may thus be written as,
\begin{equation}
  I^h=\int_{[0,L]}\left[W^h_b+W^h_e-\left( \eta^h \bs{n}^h(\tb{e}_1)-\frac{d \gamma^i_h}{d\tau}(\bar{\bs{\Lambda}}^h)^j_i n^h_j\right)-\bs{q}\cdot \bs{\gamma}^h  \right]d\tau -(\bs{\gamma}^h\cdot \bs{n})|_{\partial [0,L]}-\langle\widehat{\bs{m}},\widehat{\bs{\Omega}}^h\rangle|_{\partial [0,L]},
  \label{eq:discreteHWvartiationalPrinciple}
\end{equation}
In the above equation, $W^h_b$ and $W^h_e$ are the discrete approximations for bending and extensional energies,   $\tb{e}_1^h=\bar{\bs{\Lambda}}^h\bar{\tb{E}}_1$ and $\gamma^i_h$ denotes the
finite dimensional approximation for the components of the position vector describing the deformed configuration. An
important difference between the functional given in
\eqref{eq:HWvartiationalPrinciple} and \eqref{eq:discreteHWvartiationalPrinciple} is
that the latter need to be extremized over a finite dimensional manifold (configuration space).

\subsection{Discrete configuration space}
We introduce the discrete configuration space over which the discrete variational principle must be extremized. Even though the discrete configuration space is finite dimensional, its dimension will in general be large depending on the mesh used. As discussed in the previous subsection, for each finite element (1-simplex), we need to find approximations for frame, displacement, force 1-form and deformation 1-form. Note that the present configuration space includes not only kinematic variables but also kinetic variables as states. The inclusion of kinetic variables in the description of the configuration space is motivated by the mixed nature of our  variational principle. The discrete configuration space may thus be described as,
\begin{eqnarray}
  \mathcal{C}=\{SO(3)^n \times \left(\mathbb{R}^3\right)^n\times \left(\mathbb{R}^3 \right)^m\times \left(\mathbb{R}\right)^m\},
\end{eqnarray}
where $n$ and $m$ denote the number of nodes (0-simplices) and elements
(1-simplices) in the FE mesh. We have written down the discrete configuration
space as the product of four components; the first and second
represent the frame and displacement vectors at each node. For an FE mesh with 
 $n$ nodes, $n$ copies of $SO(3)$ and $\mathbb{R}^3$ are
required to represent the states of displacement and frame at each node. Similarly, for each finite element, the DOF's
associated with force and deformation 1-forms can be identified with
$\mathbb{R}^3$ and $\mathbb{R}$ respectively. To represent these kinetic states,
we need $m$ copies of $\mathbb{R}^3$ and $\mathbb{R}$, were $m$ is the number of elements in the FE mesh. The discrete
configuration space thus contains all the possible states that the kinetic and
kinematic variables can take.

\subsection{Discrete bending and extension energies}
Discrete approximations to the differentials of the frame field, in terms of the approximated frame fields, can be written as,
\begin{equation}
  \frac{d \bs{\Lambda}^h}{ds}=\frac{\widehat{\bs{\psi}}}{L^h}\exp\left(\frac{s}{L^h}\widehat{\bs{\psi}}\right)\bs{\Lambda}_1.
  \label{eq:diffFrame}
\end{equation}
The above relation is obtained by differentiating \eqref{eq:expInterpolation}. 
Using \eqref{eq:diffFrame}, the approximate bending strain may be expressed as,
\begin{align}
  \nonumber
  \widehat{\bs{\Omega}}^h&=(\bs{\Lambda}^h)^t\frac{d \bs{\Lambda}^h}{ds},\\
  &=\bs{\Lambda}^t_1\exp\left(-\frac{s}{L}\widehat{\bs{\psi}}\right)\frac{\widehat{\bs{\psi}}}{L}\exp\left(\frac{s}{L}\widehat{\bs{\psi}}\right)\bs{\Lambda}_1.
\end{align}
Finally, the bending energy density can be approximated as,
\begin{align}
  \nonumber
  W_b^h&=\langle \widehat{\bs{\Omega}}^h, \mathbb{C}:\widehat{\bs{\Omega}}^h \rangle_{so(3)},\\
  &=\langle \bs{\Omega}^h,\tb{D}\bs{\Omega}^h \rangle_{\mathbb{R}^3}.
  \label{eq:discretBendingEnergy}
\end{align}
The above equation is arrived at by pulling  $\langle.,.\rangle_{so(3)}$
back to $\mathbb{R}^3$ using the $\widehat{(.)}$ map. The constitutive matrix
$D$ is completely determined by $\mathbb{C}$. In general, $D$ is symmetric and
positive definite. Similarly, the discrete approximation for extensional energy can be written as,
\begin{eqnarray}
  W_e^h=\frac{AE}{2}(\eta^h-1)^2.
\end{eqnarray}

\subsection{Residue and tangent}
The critical points of the discrete functional given in \eqref{eq:discreteHWvartiationalPrinciple} can be numerically computed through any optimization algorithm. In this work, we use Newton's method, which is a second order optimization technique. Newton's method seeks to improve the solution by producing a sequence of increments. Computing each increment requires the evaluation of both first and second
derivatives of the objective functional which are commonly known as the residue
vector and the tangent operator.  Our choice of Newton's method as the optimizer stems from our desire to track the equilibrium path for given forcing conditions. The discrete equilibrium in terms of the differential of the energy functional can be written as,
\begin{eqnarray}
  \ed I^h=0.
  \label{eq:residue}
\end{eqnarray}
In the finite element literature, $\ed I^h$ is often called the residue. Following this convention, we may also denote it by $\mathcal{R}$. In the present context, $\ed I^h$ is a 1-form from $T^*\mathcal{C}$. To calculate the components of the residue, we use the Gateaux differential. The residue is constructed by stacking Gateaux differentials of $I^h$ in multiple directions,
\begin{equation}
  \mathcal{R}:=[\ms{D}_{\bs \eta} I^h, \;
  \ms{D}_{\tb n^1}      I^h, \; 
  \ms{D}_{\tb n^2}      I^h, \; 
  \ms{D}_{\tb n^3}      I^h, \; 
  \ms{D}_{\tb u^1} I^h,\;
  \ms{D}_{\tb u^2} I^h,\;
  \ms{D}_{\tb u^3} I^h,\;
  \ms{D}_{\bs \lambda^1} I^h,\;
  \ms{D}_{\bs \lambda^2} I^h,\;
  \ms{D}_{\bs \lambda^3} I^h]^t.
  \label{eq:residue}
\end{equation}
In the last equation, $\bs{\eta}$, $\tb{n}^i$, $\tb{u}^i$ and $\bs{\lambda}^i$ denote the DOF's associated with deformation 1-forms, force 1-forms, displacements and rotation fields respectively.   The Gateaux derivative of $I^h$ in the direction of a DOF is denoted by $\ms D_{(.)}I^h$; expressions for these components can be found in Appendix \ref{sec:residueTangent}. The second derivative of $I^h$ should take into account the nonlinear nature of the configuration manifold $\mathcal{C}$. Simo \cite{simo1992symmetric} proposed a technique to calculate the correct (symmetric) Hessian for an energy function when the discrete configuration space is endowed with a metric. Incorporating this technique in our calculations, the tangent stiffness can be written as,
\begin{equation}
  \mathcal{K}=
  \begin{bmatrix}
      \ms D_{\bs \eta\bs \eta}I^h & \ms D_{\bs \eta\tb n^j}I^h &\ms D_{\bs \eta\tb u ^j}I^h &\ms D_{\bs \eta\bs \lambda^j}I^h\\
       &\ms D_{\tb n^i\tb n^j}I^h & \ms D_{\tb n^i\tb u^j}I^h  & \ms D_{\tb n^i\bs \lambda^j}I^h\\ 
      & & \ms D_{\tb u^i\tb u^j}I^h & \ms D_{\tb u^i\bs \lambda^j}I^h\\
       \text{Symm}&  &  &  \ms D_{\bs \lambda^j\bs \lambda^j}I^h
  \end{bmatrix};\quad i,j=1,...,3.
  \label{eq:tangent}
\end{equation}
Expressions for the components of the stiffness matrix $\ms D_{(.)(.)I^h}$ may be found in Appendix \ref{sec:residueTangent}. 

\subsection{Configuration update}
Solution update is an important aspect that needs care when the configuration space is not an Euclidean space. An iterative optimization procedure (Gradient descent, Newton's
method, etc.)  produces an incremental update to the last known solution of the extremization problem. This increment can be identified with the tangent space of the configuration manifold. If the configuration space were Euclidean, just adding the solution to the incremental solution would produce a valid configuration. In the present work, we have the frame as one of our variables whose evolution is defined on $SO(3)$. Since adding an element from the tangent space of $SO(3)$ does not produce an element in $SO(3)$, the conventional additive update will not serve our purpose. The additive update has to be replaced by a more appropriate scheme, e.g. the retraction map \cite{absil2009optimization}. The idea of retraction was introduced in the literature on constrained optimization to keep the incremental solution consistent with the constraints imposed in the search space of the optimization problem. The exponential map on $SO(3)$ is an example of retraction that sends the tangent space to the manifold. Often  finding explicit expressions for a retraction is challenging; however in the present case we have an explicit expression for the exponential map in $SO(3)$.
Simo and Vu-Quoc \cite{simo1986three} introduced the exponential update on $SO(3)$ to increment the solution within a Newton-Raphson solution scheme. In this work, we use the exponential update scheme proposed by Simo and Vu-Quoc \cite{simo1986three}. An algorithm describing Newton's method along with the exponential update scheme is shown below.

\framebox[0.95\width]{

\begin{algorithm}[H]
\begin{center}
\textbf{Algorithm:} Newton's method for extremizing the discrete mixed functional
\end{center}
	\SetAlgoLined
	 \KwData{\\
      1. Finite element approximation for $ \Gamma$\;
      2. Deformation 1-form $\bs{\eta}^1_{prev}$ from the previous load-step\; 
      3. Traction 1-form $\tb{n}^i_{prev}$ from the previous load-step\;
      4. Displacement $\tb{u}^i_{prev}$ from the previous load-step\;
      5. Rotation $\bs{\Lambda}^i_{prev}$ from the previous load-step\;
      6. Displacement and traction boundary conditions for the current load step\;
	}
	 \KwResult{$\bs{\eta}^1, \tb{t}^i, \tb{u}^i$ and $ \bs{\lambda}^i$ for the prescribed displacement and traction conditions}
	 Initialize  $\mathcal{K}$ as a sparse matrix and $\mathcal{R}$ as a sparse vector\;
     Initialize  $TOL$\;
     $\tb{u}^i\leftarrow\tb{u}^i_{prev}$\;
     $\tb{n}^i\leftarrow\tb{n}^i_{prev}$\;
     $\bs{\eta}^1\leftarrow\bs{\eta}^1_{prev}$\;
     $\bs{\Lambda}^i\leftarrow\bs{\Lambda}^i_{prev}$\;

    Modify $\mathcal{R}$ to incorporate traction and displacement conditions\;

	 \While{$\norm{\mathcal{R}} \geq TOL$}{
		\tb{Initialization:} $itEl=1$\\	
		
		 \While{itEl $\leq$ \# Elements}{
                Compute element $R$ using  \eqref{eq:residue}\;
                Compute element $K$\ using \eqref{eq:tangent}\;
                Assemble $\mathcal{R}:\mathcal{R}\leftarrow R$\;
                Assemble $\mathcal{K}:\mathcal{K}\leftarrow K$\;
                }
                Apply displacement boundary conditions to $\mathcal{K}$\;
                Solve: $\mathcal{K}X=-\mathcal{R}$\;
                $\bs{\eta}^1\leftarrow\bs{\eta}^1+\Delta\bs{\eta}^1$\;
                $\tb{n}^i\leftarrow\tb{n}^i+\Delta\tb{n}^i$\;
                $\tb{u}^i\leftarrow\bs{u}^i+\Delta\bs{u}^i$\;
                $\bs{\Lambda}^i\leftarrow \exp({\bs{\lambda}^i})\bs{\Lambda}^i$\;
	  
	 }
         \label{algo:NewtonsMethod}
	\end{algorithm} 
}\par



\section{Numerical results}
\label{sec:NumericalResults}

In this section, we apply the numerical technique discussed in Section \ref{section6:Finite element approximation} to large deformation simulations of Kirchhoff  rods. The problems presented in this section are chosen to test the robustness of our   numerical approximation. First, we test the frame indifference of our implementation by rigidly rotating an L-shaped cantilever beam. The rate of convergence of our numerical approximation is then assessed using a cantilever subjected to end shear and moment. Williams frame is studied for its snap through behaviour, while the circular arch is studied for its complicated load deformation behaviour. The bent cantilever is studied to  demonstrate our ability to solve a curved reference configuration undergoing large deformation. Similarly, the L-shaped cantilever is studied for non-smooth reference configuration (reference configuration with corners). Our ability to solve problems with extremely large rotations is demonstrated using a cantilever subjected to end torque and moment. Finally, our ability to solve structural systems with multiple member intersections is demonstrated with a star shaped dome.  

In the numerical results reported below, the bending energy is assumed to be quadratic in curvature  leading to a linear curvature bending moment relationship. The integrals involved in the computation are performed using a two point Gaussian quadrature rule. Newton's method is used to extremize the discrete mixed variational functional. The Newton iteration is terminated when the norm of the residue becomes smaller than $10^{-5}$. 

\subsection{Frame indifference}
\begin{figure}
  \begin{subfigure}{.9\textwidth}
    \centering
    \includegraphics[scale=0.8]{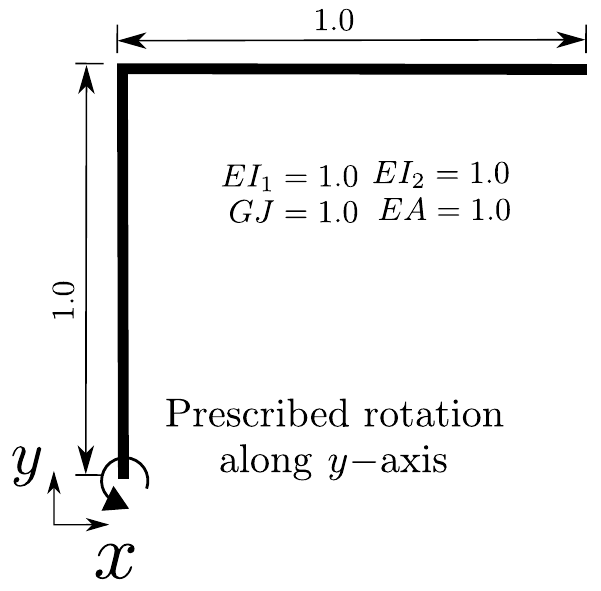}
    \caption{}
 \label{fig:frameIndiffBvp}
  \end{subfigure}%
  \\
  \begin{subfigure}{.5\textwidth}
    \centering
    \includegraphics[scale=1.3]{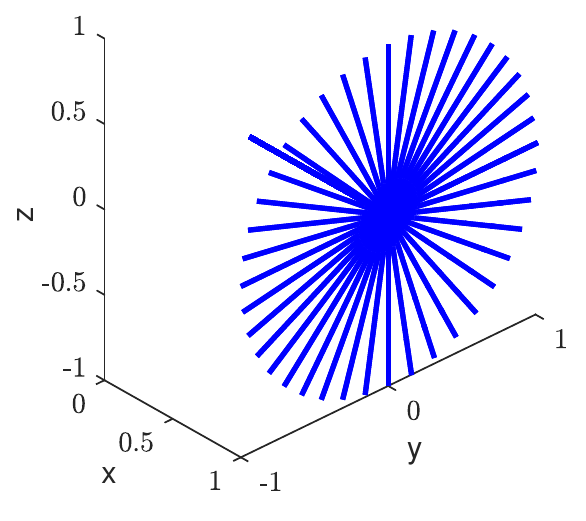}
    \caption{}
 \label{fig:frameIndiff1}
  \end{subfigure}%
  \begin{subfigure}{.5\textwidth}
    \centering
    \includegraphics[scale=0.8]{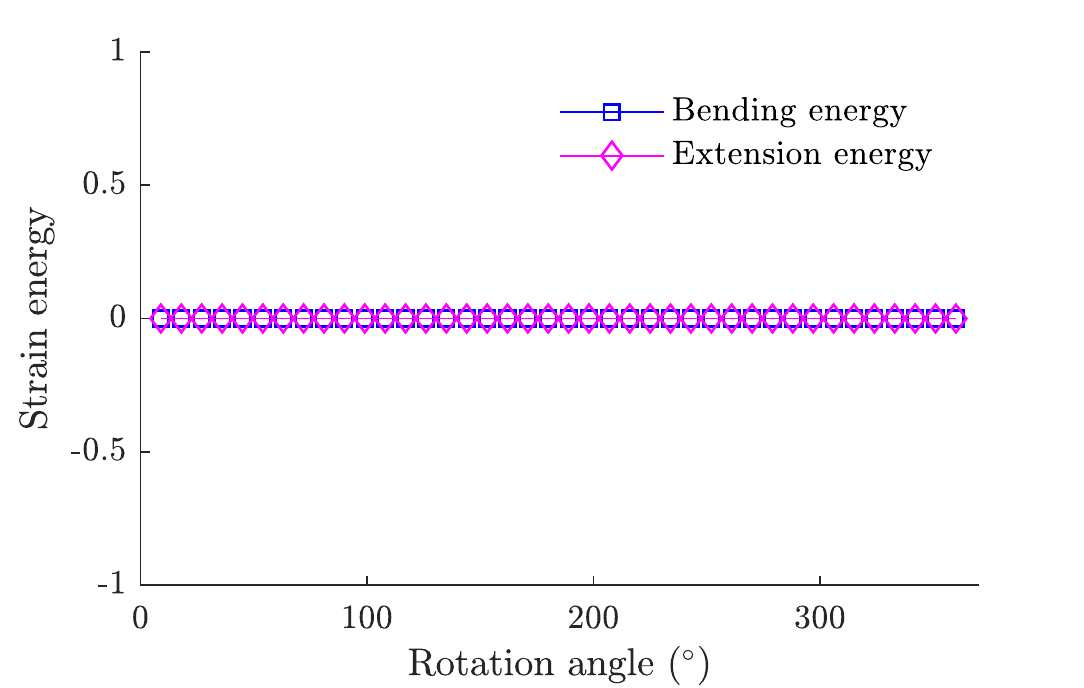}
    \caption{}
 \label{fig:frameIndiff2}
  \end{subfigure}
  \caption{(a) L-frame with prescribed rotation boundary condition; (b) sequence of deformed configurations for prescribed end rotation; (c) variation of strain energy with respect to angle of rotation.}
  \label{fig:frameIndiff}
  \end{figure}

We now demonstrate that the present numerical implementation is frame indifferent.
Towards this we study the rigid rotation of an L-frame about the $x$ axis. Displacement DOF at one end of the frame is fixed while the other end is left free. A non-zero rotation boundary condition is applied at the same end where the zero displacement boundary condition is applied.  Since the L-frame experiences no externally
applied loads, its deformed configuration is a rigid rotation about the $x-$axis.  Additional Lagrange multipliers are introduced to impose the rotation
constraint at the fixed end. The modified mixed functional accounting for the rotation constraint can be written as,
\begin{equation}
  I^h=\int_{[0,L]}\left[W^h_b+W^h_e-\left( \eta^h \bs{n}^h(\tb{e}_1)-\frac{d \gamma^i_h}{d\tau}(\bar{\bs{\Lambda}}^h)^j_i n^h_j\right)-\bs{q}\cdot \bs{\gamma}^h  \right]d\tau -(\bs{\gamma}^h\cdot \bs{n})|_{\partial [0,L]}-\langle\widehat{\bs{m}},\widehat{\bs{\Omega}}^h\rangle|_{\partial [0,L]},+g(\bs{\Lambda}_n,\bs{\Lambda},\bs{\rho}),
\end{equation}
where, $g=\bs{\rho}:(\bs{\Lambda}^t_n\bs{\Lambda}_c-\bs{I})$ and $:$ denotes the
inner-product between two matrices, $\bs{\Lambda}_n$ is the rotation prescribed at the fixed end, $\bs{\Lambda}_c$ is the computed rotation and $\bs{\rho}$ is the skew-symmetric matrix introduced as a Lagrange multiplier to impose the rotation constraint. The components of residue and tangent operators can be computed using the same linearization procedure presented in the previous section. The deformed configurations
predicted by the present numerical procedure for different values of prescribed
rotation is presented in Figure \ref{fig:frameIndiff1}. Bending and extension
energies stored in the beam are shown in Figure
\ref{fig:frameIndiff2}; both of which are zero. This exercise clearly demonstrates that, along with the beam formulation, the numerical implementation is frame indifferent too.   

\subsection{Rolling and unrolling of a cantilever beam}
\begin{figure}
  \begin{subfigure}{.5\textwidth}
    \centering
    \includegraphics[scale=0.3]{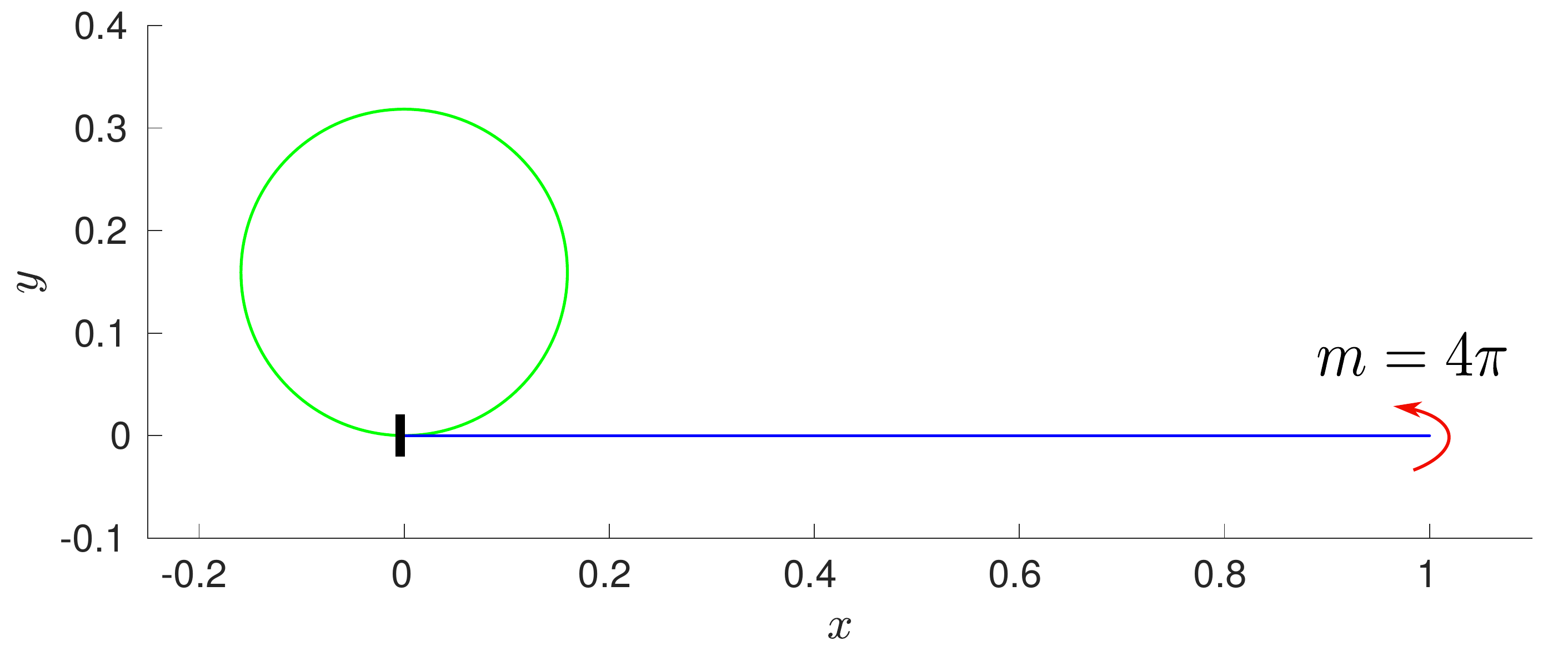}
    \caption{}
    \label{fig:rollingBeam1}
  \end{subfigure}%
  \begin{subfigure}{.5\textwidth}
    \centering
    \includegraphics[scale=0.8]{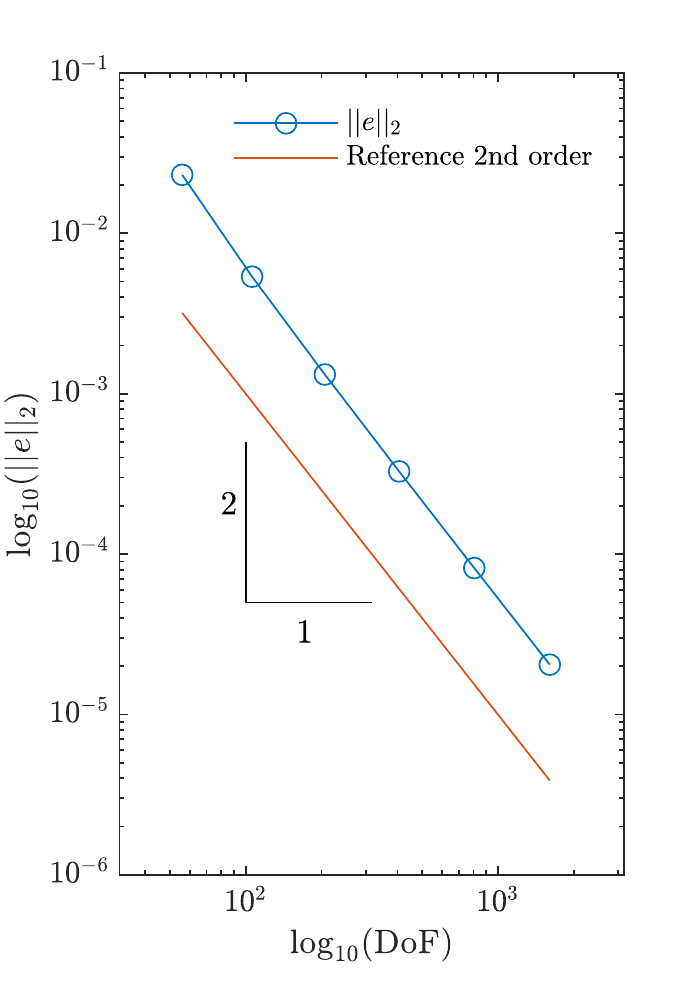}
    \caption{}
 \label{fig:rollingBeam2}
  \end{subfigure}
  \caption{(a) Deformed configuration of the cantilever subjected to end moment, the blue curve indicates the reference configuration, while the green curve indicates the deformed configuration. (b) Convergence of deformation to the analytical solution.}
  \label{fig:rollingBeam}
  \end{figure}

\begin{figure}
    \begin{subfigure}{.45\textwidth}
      \centering
      \includegraphics[scale=0.35]{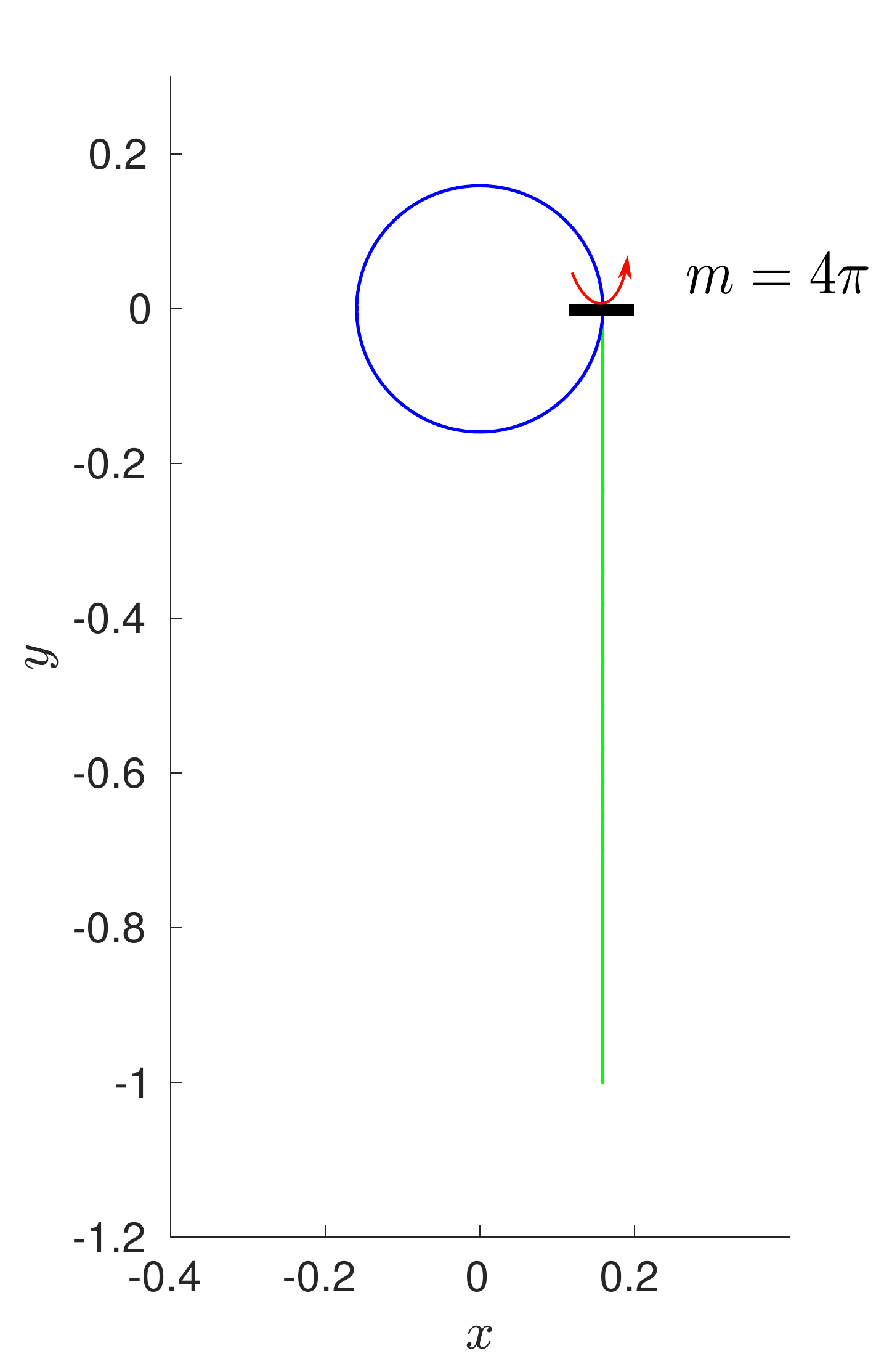}
      \caption{}
\label{fig:unRollingBeam1}
    \end{subfigure}%
    \begin{subfigure}{.45\textwidth}
      \centering
      \includegraphics[scale=0.8]{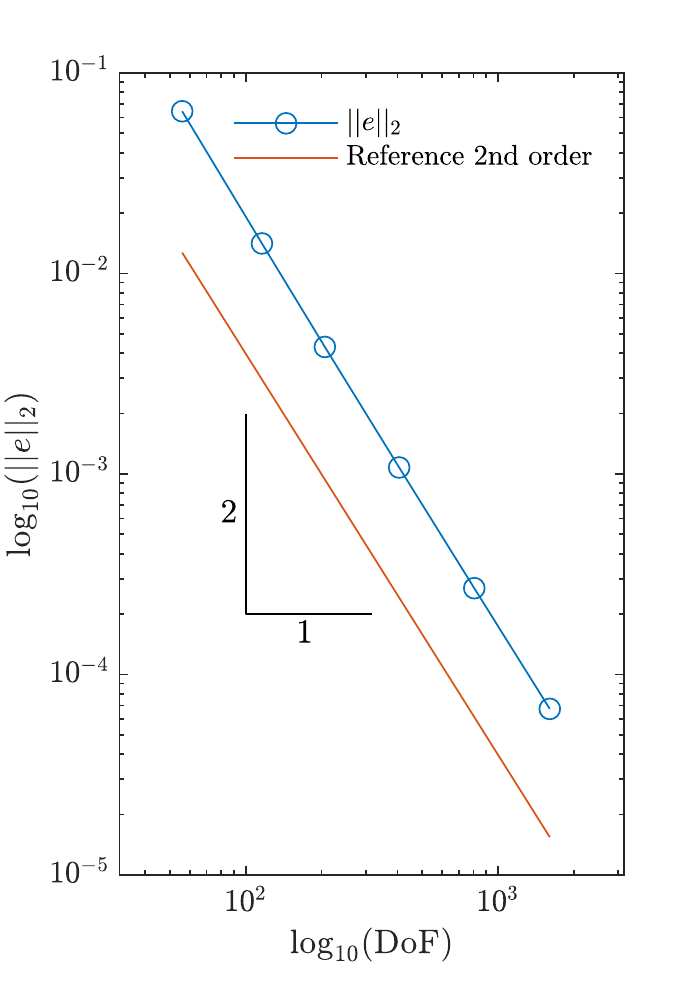}
      \caption{}
\label{fig:unRollingBeam2}
    \end{subfigure}
    \caption{(a) Deformed configuration of the unrolled beam; the blue curve indicates the reference configuration and the green curve the deformed configuration. (b) Convergence of deformation to the analytical solution}
    \label{fig:unRollingBeam}
    \end{figure}

We now study the performance of our numerical scheme against pure bending problems. We consider two cases; while the first case studies the rolling of a straight cantilever beam, the second case involves unrolling a curved cantilever beam.  For both the cases, axial stiffness of the beam is chosen as $AE=10$ and bending and torsional stiffnesses $(EI_1, EI_2$ and $GJ)$ are taken as 2. To evaluate the convergence of the FE approximation, we use the relative $L_2$ error. Discretizations with 5, 10, 20, 40, 80 and 160 finite elements were used to study the convergence. The $L_2$ error in deformation is given by,
\begin{equation}
{\| e \|}_2=\sqrt{ \int_{0}^{L} {{\| \bs{\gamma}-\bs{\gamma}^{h} \|}^2} \,d\tau },
  \label{eq:norm}
\end{equation}
where $\bs{\gamma}$ and $\bs{\gamma}^{h}$ are the centerlines of the deformed rod via the original formulation and its discrete (numerical) representation. Figure \ref{fig:rollingBeam1} shows the deformed configuration predicted by our numerical implementation for a straight cantilever beam subject to an applied end moment of $4\pi$. In Figure \ref{fig:rollingBeam2}, we present the FE convergence of the rod's centerline. From the convergence plots, it can be seen that the convergence is second order, which is also the optimal order for an elliptic second order system discretized using linear finite elements. Now, we study the convergence of the FE approximation for the case of a beam with curved reference configuration. A cantilever beam bent in the form of a circle with radius $1/2\pi$ is taken as the initial configuration. The same bending moment as in the straight case was applied to the free end of the beam. Figure \ref{fig:unRollingBeam1} shows the reference and deformed configurations along with the boundary conditions. Because of the applied bending moment, the circular beam unrolls to become straight.  Figure \ref{fig:unRollingBeam2} show the convergence of deformation. As in the rolling problem, the order of convergence for the unrolling problem was found to be two.

\subsection{Cantilever beam subjected to shear load}
\begin{figure}
  \begin{subfigure}{.45\textwidth}
    \centering
    \includegraphics[scale=0.450]{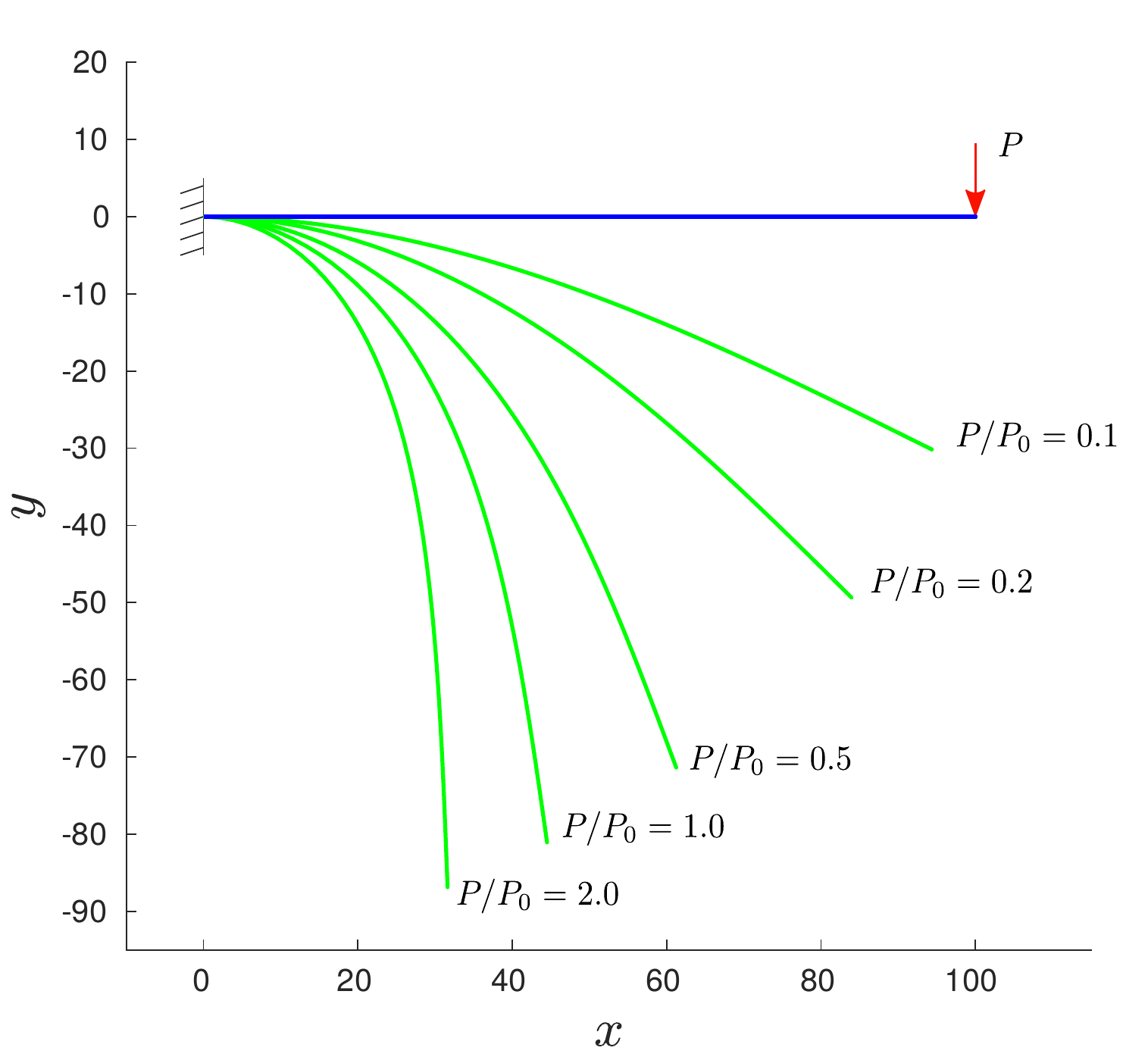}
    \caption{}
  \label{fig:elastica1}
  \end{subfigure}%
  \begin{subfigure}{.45\textwidth}
    \centering
    \includegraphics[scale=0.31]{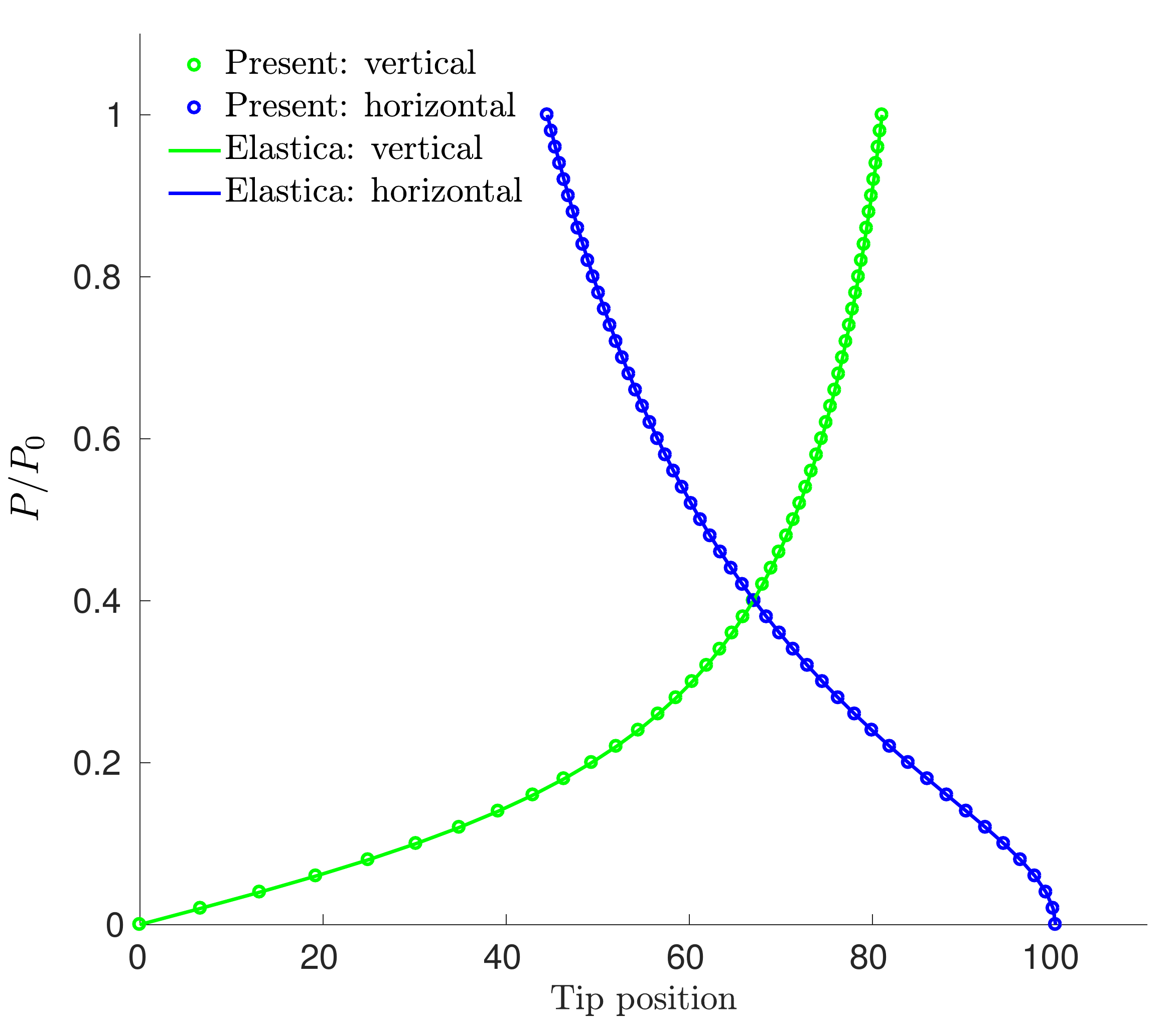}
    \caption{}
  \label{fig:elastica2}
  \end{subfigure}
  \caption{(a) Deformed configuration of the cantilever beam predicted by the present mixed FE model; (b) load displacement plot for vertical and horizontal displacement}
  \label{fig:elastica}
  \end{figure}

 We now demonstrate the performance of the present FE formulation under shear
 load. As a benchmark, we study the elastica, whose analytical solution can be expressed in terms of elliptic integrals. These solutions may be found in
 \cite{frisch1962flexible}.  In the simulation, we chose the
 length of the rod to be 100 and $EI$ to be 1000.  The beam undergoes large
 displacements due to the applied shear. In Figure \ref{fig:elastica1}, the
 reference configuration is plotted along with the simulated deformed configurations via our
 mixed FE model for various magnitudes of applied shear. A
 comparison of the numerically computed tip displacement and elastica
 solution is presented in Figure \ref{fig:elastica2}.  
 From the figure, one observes that  our FE formulation accurately predicts the
 tip displacement.
 
\subsection{Williams frame}
\begin{figure}
\begin{subfigure}{.45\textwidth}
  \centering
  \includegraphics[scale=0.850]{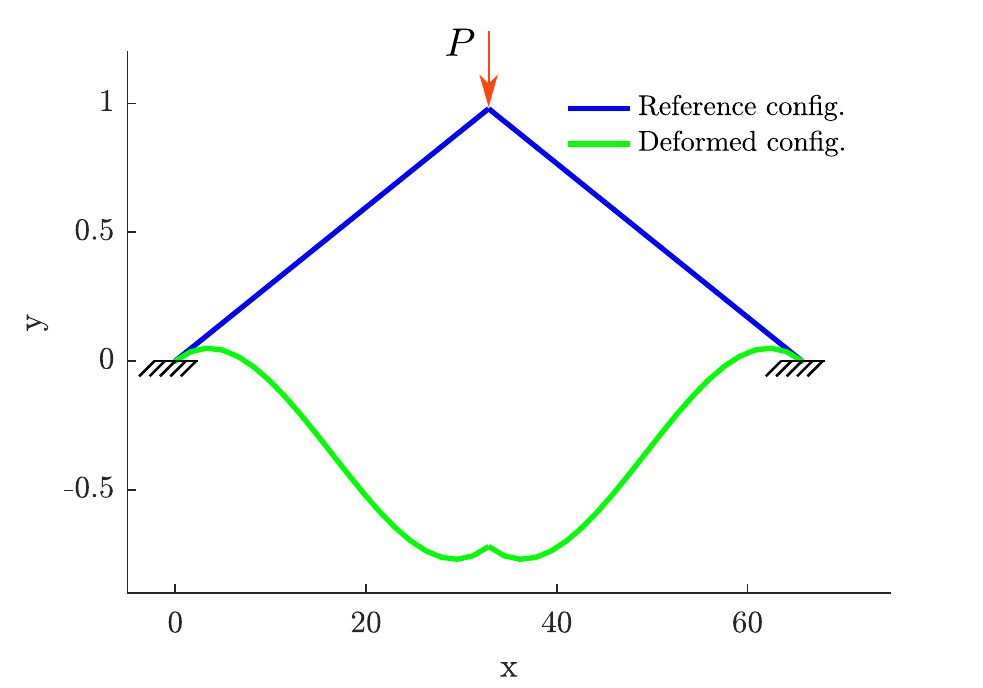}
  \caption{}
\label{fig:toggleFrame1}
\end{subfigure}%
\begin{subfigure}{.45\textwidth}
  \centering
  \includegraphics[scale=0.850]{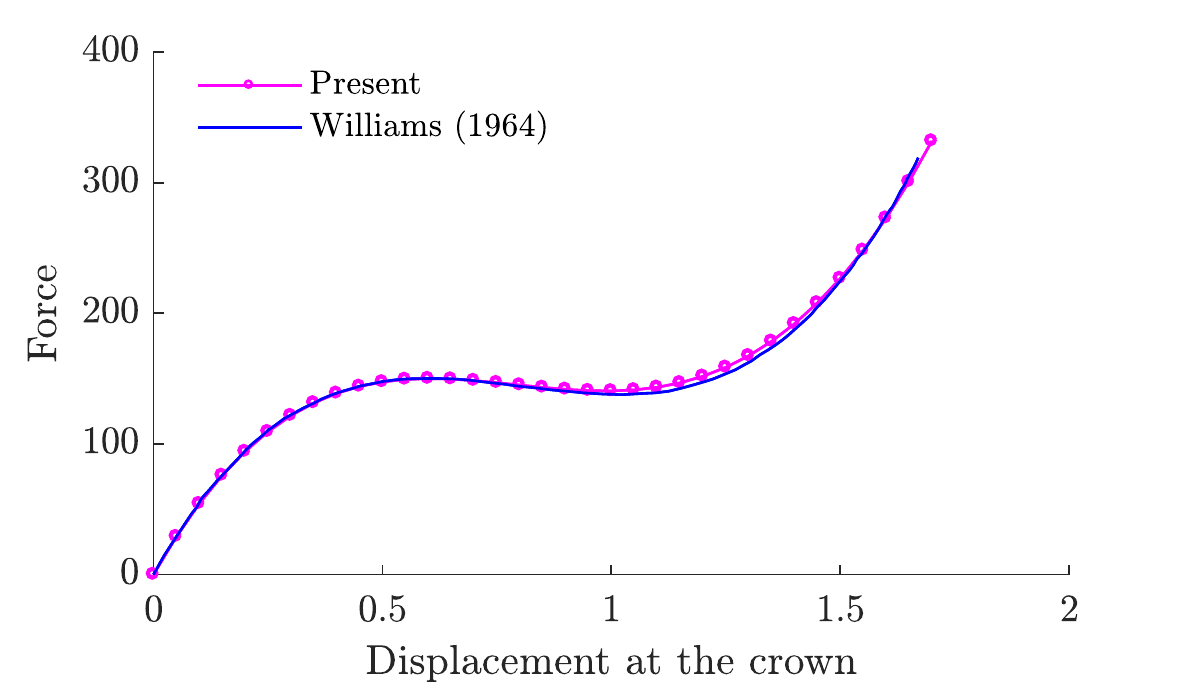}
  \caption{}
\label{fig:toggleFrame2}
\end{subfigure}
\caption{(a) Deformed configuration of the Williams frame predicted by the present FE formulation; (b) force versus crown displacement predicted using a displacement controlled loading.}
\label{fig:toggleFrame}
\end{figure}

Williams frame is a classical post buckling problem \cite{williams1964approach},  studied using a Kirchhoff  rod theory. Considered as a benchmark problem to assess numerical techniques for non-linear beam formulations, this problem has attracted widespread attention in the literature.  
The frame contains two inclined beam members connected rigidly at the crown,
while the other ends of both the beams are fixed rigidly, i.e) both displacement and rotation are constrained. A vertical force in the
plane of the frame is applied at the crown. Figure  \ref{fig:toggleFrame1} shows the reference configuration along with loading and boundary conditions. The span of the toggle frame between the supports is taken as 65.715 units, while the crown is located at a height of 0.98 units from the support. The limbs of the frame are assumed to have a circular cross section with diameter 0.721 units and Young's modulus 199714 units. The load displacement curve at the crown predicted by our mixed FE model is shown in Figure \ref{fig:toggleFrame2}. A displacement controlled loading is used to the trace the load displacement curve presented in the figure. 
Under the action of the applied load, the crown deflects vertically. After reaching
a critical load, the frame toggles, causing a slight reduction in load. Our displacement controlled solution technique captures the unstable path of the load displacement plot
as well. In Figure \ref{fig:toggleFrame2}, the load displacement predicted by
\cite{williams1964approach} is plotted alongside the present predictions for a
ready comparison and a good agreement can be found.      

\subsection{Circular arch}
\begin{figure} [h]
  \begin{subfigure}{.45\textwidth}
    \centering
    \includegraphics[scale=0.30]{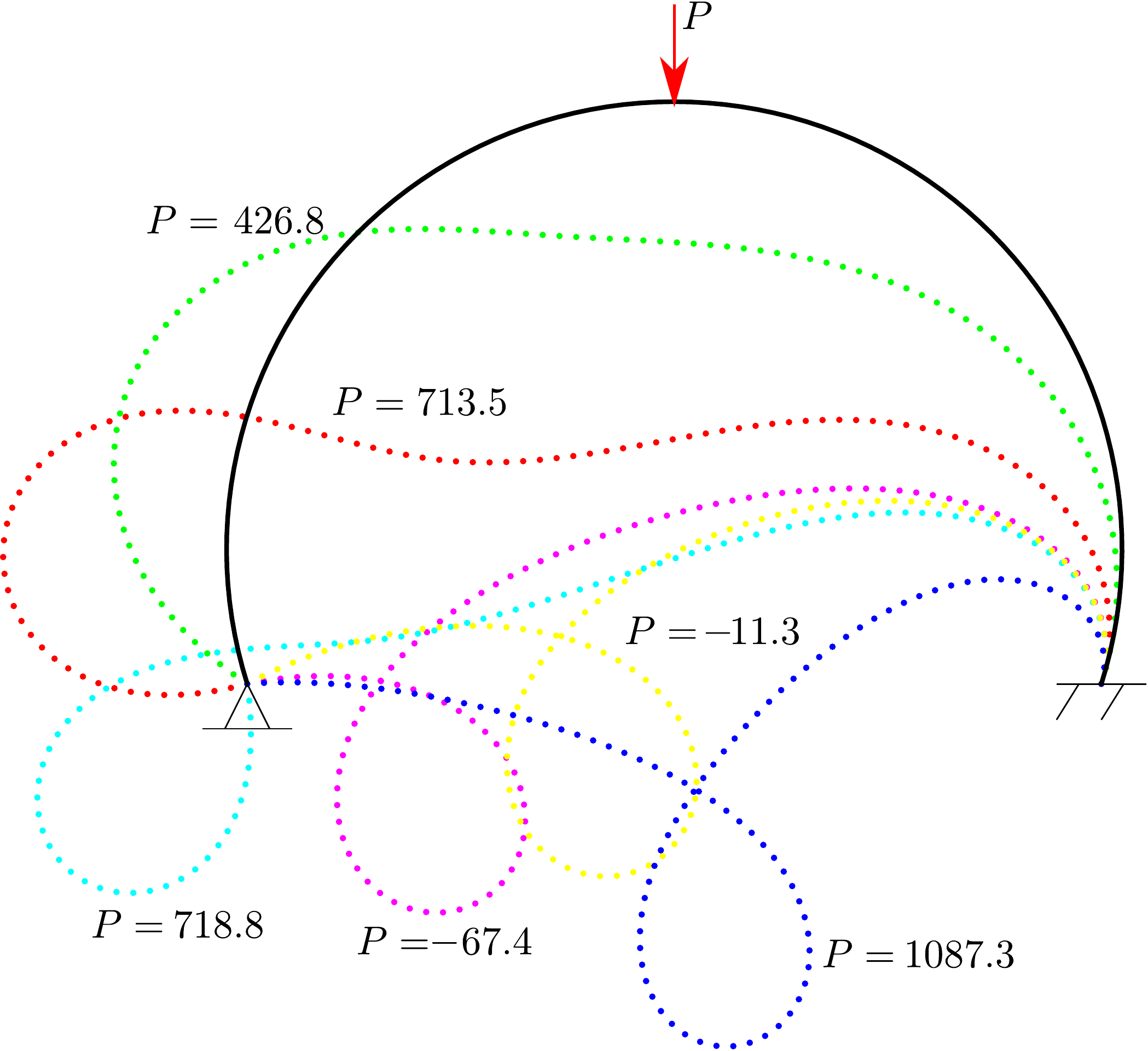}
    \caption{}
    \label{fig:circArch1}
  \end{subfigure}%
  \begin{subfigure}{.45\textwidth}
    \centering
    \includegraphics[scale=0.25]{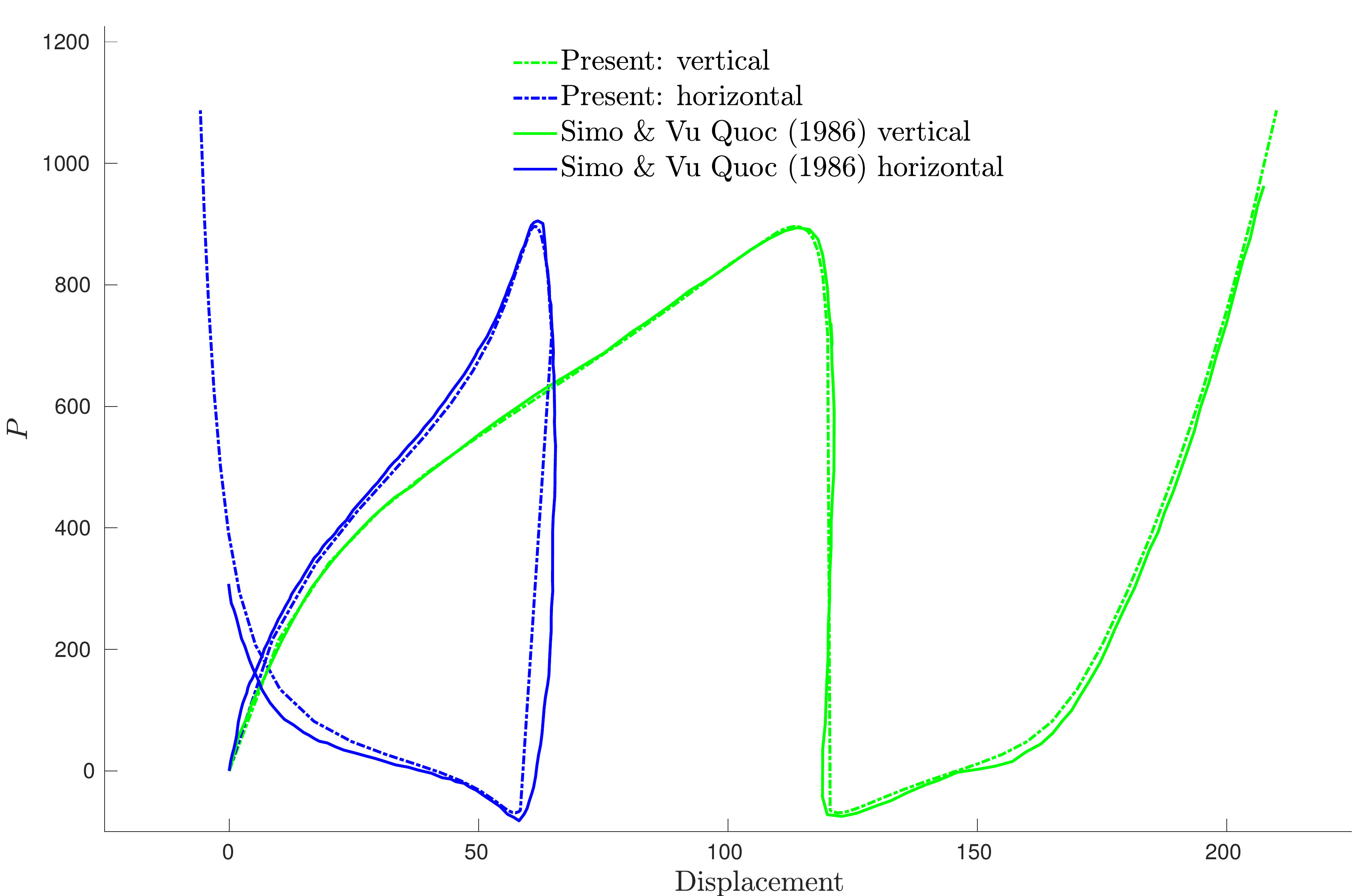}
    \caption{}
    \label{fig:circArch2}
  \end{subfigure}
  \caption{(a) Arch subjected to a concentrated load at the crown:  sequence of deformed configurations; (b) load displacement comparison for vertical and horizontal displacements}
  \label{fig:circArch}
  \end{figure}

We now study the present formulation fares for post bucking
problems where the reference configuration is curved. Thus consider a deep circular
arch subjected to a shear load at the crown. This problem has been studied by may
authors. Simo and Vu-Quoc \cite{simo1986three} first traced the complete post-bucking path of the deep ciruclar arch using a Reissner rod. Notably, this problem has multiple load
paths emanating from a bifurcation point; we only track the load path presented in Simo and Vu-Quoc. 
The arch has a radius of 100 units and interior angle of 215$\deg$. The
constants $EI_1$, $EI_2$ and $GJ$  are set as $10^6$, while $AE$ is taken to be $10^8$. The arch is clamped at one end and pinned at the other with a  central point load. Figure
\ref{fig:circArch1} shows the reference configuration along with the boundary conditions. The
initial configuration is discretized using 100 finite elements. The numerical
simulation was performed using a displacement controlled strategy. A small horizontal load of 8 units (1\% of the first buckling load)  is
applied at the crown to force the solution procedure to choose the present load
path. Figure \ref{fig:circArch1} shows the sequence of deformed configurations
predicted by our numerical technique for different vertical loads.


The load displacement (at crown) predicted by our numerical technique is
plotted in Figure \ref{fig:circArch2}, which shows both vertical and horizontal
displacements. The present analysis predicts a
critical buckling load of 896.13, while \cite{dadeppo1975instability} predicts it
as 897.   DaDeppo and Schmidt \cite{dadeppo1975instability} based their buckling analysis on Euler's inextensible elastica. The initial buckling load estimated by
\cite{dadeppo1975instability} compares closely with the present approach since
both are based on the Kirchhoff theory. The predictions of Simo and Vu-Quoc \cite{simo1986three} are slightly higher since their analysis is based on a Reissner  beam theory. In Figure \ref{fig:circArch2}, the load deformation plot predicted by our scheme is plotted along with the prediction by \cite{simo1986three}. From the load displacement plot, we find that the
our mixed FE model captures the unstable branch of the load
displacement curve quite well. A good comparison is obtained for the final
stable branch as well.

\subsection{Bent cantilever}
\begin{figure}
  \begin{subfigure}{.45\textwidth}
    \centering
    \includegraphics[scale=0.40]{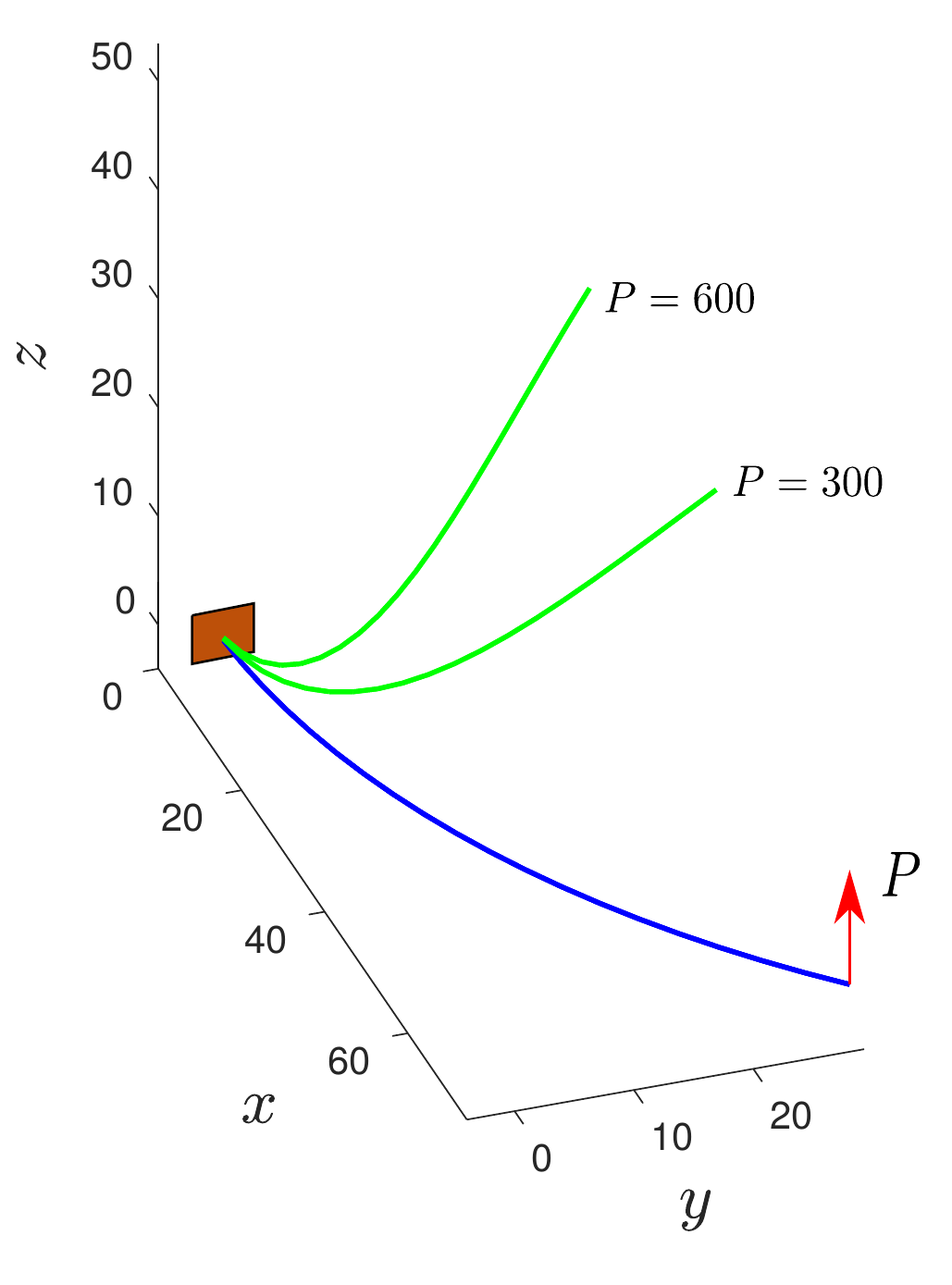}
    \caption{}
    \label{fig:bentCantilever1}
  \end{subfigure}%
  \begin{subfigure}{.45\textwidth}
    \centering
    \includegraphics[scale=0.8]{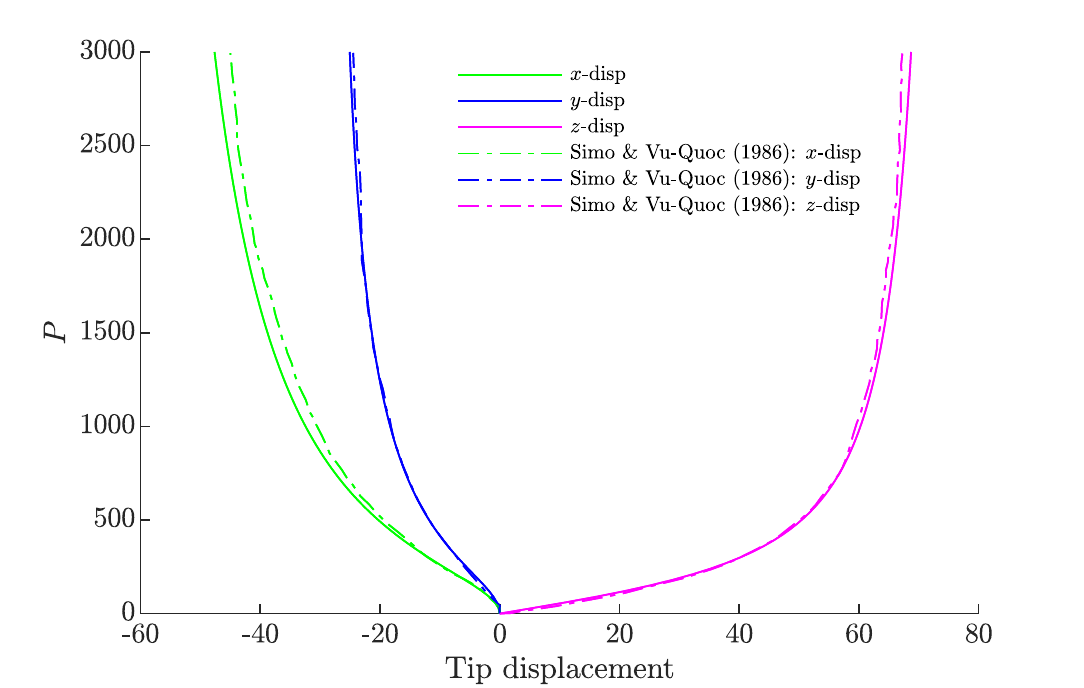}
    \caption{}
    \label{fig:bentCantilever2}
  \end{subfigure}
  \caption{(a) Deformed configuration of the bent cantilever for different load levels;  (b)  comparison of load displacement with Simo and Vu-Quoc (1986)}
  \label{fig:bentCantilever}
  \end{figure}

We now apply our mixed FE model to the case of a rod whose deformation is not contained to a plane.  A cantilever beam with initial configuration in the
form of a quarter circular arc is studied. Figure \ref{fig:bentCantilever1}
shows the initial configuration. The radius of the circular arc is 100. The
curved beam has one of its ends fixed, while the other end is subjected to a
shear force applied in the direction perpendicular to the plane of the circular
arc. This shear force causes the beam to undergo shear, bending and torsion
simultaneously.  The constitutive matrix relating the bending moment and
curvature is set as Diag$[10^6,10^6,10^6]$, while $AE$ is $10^3$.   A
comparison of the tip displacement predicted by our numerical implementation vis-\'a-vis data from the literature is shown in Table \ref{tab:bentCantilever} for a
shear force of 300 and 600 units. From Table \ref{tab:bentCantilever}, it is found that the tip displacement predicted by the present technique
agrees well with those reported in the literature. The load
displacement curve at the tip predicted for extremely large deformation is shown in Figure
\ref{fig:bentCantilever}. For a comparison, the load displacement curve obtained by
\cite{simo1986three} is also plotted and a good
agreement is found in these plots as well.

\begin{table}
  \centering
  \caption{Comparison of tip coordinates for different load levels}
  \begin{tabular}{l| c| c}
    \hline
    Authors & 300 & 600\\
    \hline
    Present    &  58.80, 22.26, 40.19    &  47.18, 15.70, 53.50\\
   Ghosh and Roy \cite{ghosh2009frame} &   58.89, 22.26, 40.08  &  47.29, 15.67, 53.37 \\
   Simo and Vu-Quoc \cite{simo1986three} & 58.84, 22.33, 40.08 &  47.23, 15.79, 53.37\\
   M\"akinen \cite{makinen2007total} &- &  47.01, 15.62, 53.50\\
   Mohan \cite{mohan2011group}& 58.69, 22.19, 40.33  &  47.14, 15.67, 53.48 \\
   Cardona and Geradin \cite{cardona1988beam}& 58.64, 22.14, 40.35 & 47.04,  15.55, 53.50\\
    \hline
   \end{tabular}
   \label{tab:bentCantilever}
\end{table}

\subsection{L-shape cantilever}
\begin{figure}
\centering
\begin{subfigure}{.5\textwidth}
  \centering
  \includegraphics[scale=0.5]{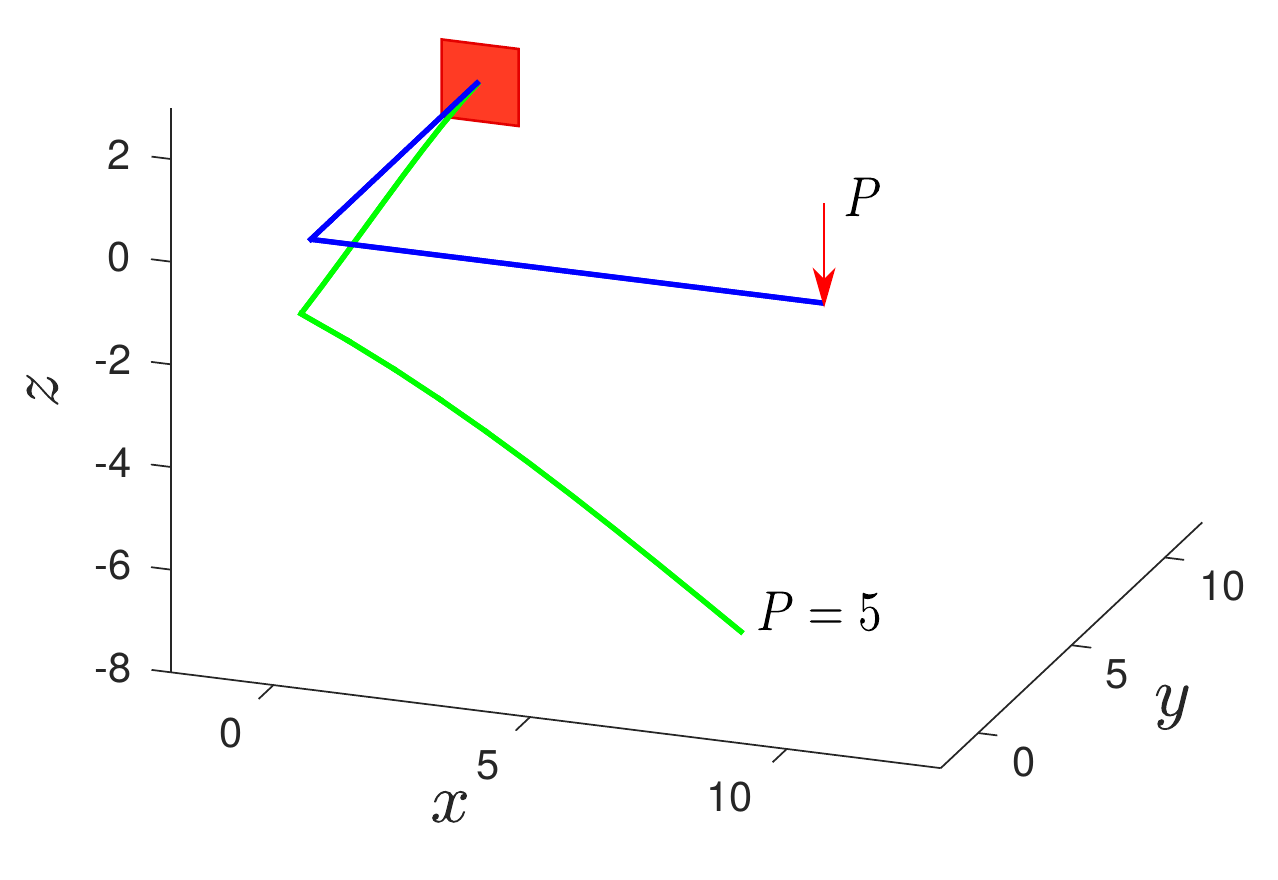}
    \caption{}
  \label{fig:frameinvariance1}
\end{subfigure}%
\begin{subfigure}{.5\textwidth}
  \includegraphics[scale=0.5]{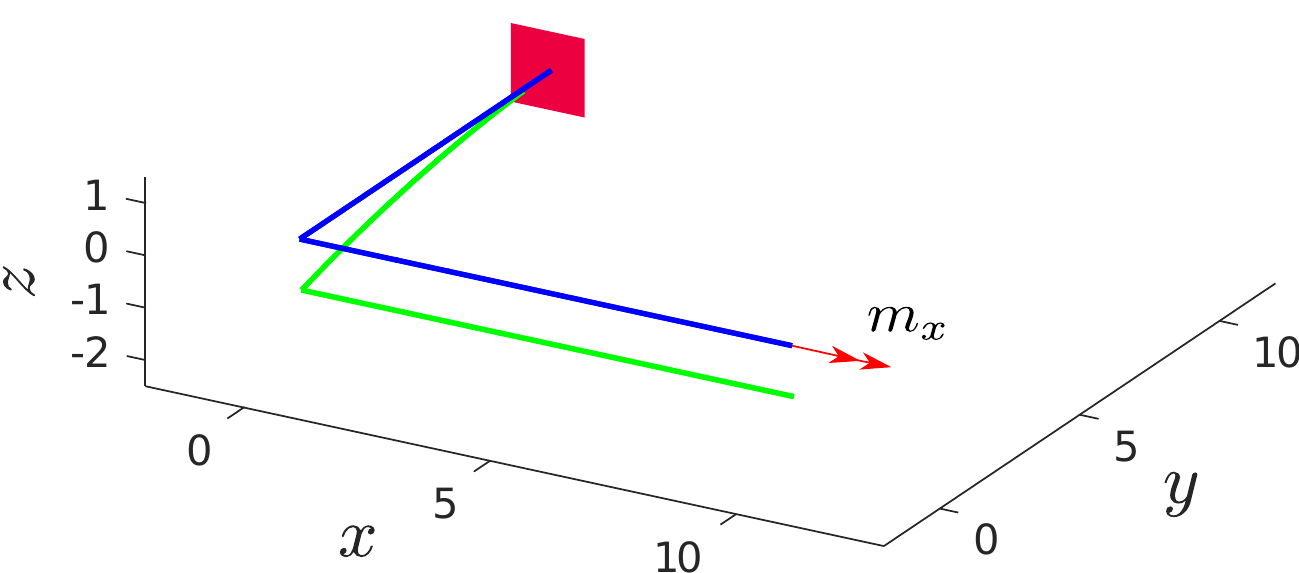}
    \caption{}
  \label{fig:frameinvariance2}
\end{subfigure}
\caption{Deflected shape of L shape cantilever under a (a) tip load and (b) torsion}
\label{fig:frameinvariance}
\end{figure}

In this subsection we consider the deformation of a frame consisting of two straight rods attached orthogonal to each other. The connection between the rods is assumed to be moment resisting. The reference configuration of the frame assembly is shown in Figure \ref{fig:frameinvariance1}. Simulations with two different boundary conditions are performed, the stiffness  properties of both the rod are taken as $AE=10^6$ and $EI_1= EI_2= GJ=10^3$. In both the simulations, the initial configuration is discretized using 8 finite elements per member with one of the free ends fixed against both displacement and rotation. The first simulation predicts the deformation of the frame due to an  applied shear force of 5 units at the free end. The deformed shape predicted by the simulation for the applied shear case is shown is Figure \ref{fig:frameinvariance1}. A comparison of the computed position vector of the tip in the deformed configuration is shown in Table \ref{tab:frameinvariance_table}, which agrees well with the the data available from the literature. The second simulation predicts the deformation of the frame subjected to a torque applied at the free end. The deformed configuration predicted by our numerical scheme is shown is Figure \ref{fig:frameinvariance1}. It should be noted that in both the simulations, no additional constraints were introduced to transfer the forces from one member to the other. The usual finite element assembly was sufficient to take care of this force transfer.
\begin{table}
  \centering
  \caption{L frame :Comparison of tip displacement data}
  \begin{tabular}{l| c|c}
    \hline
    Authors & Displacement($x,y,z$) & \# elements\\
    \hline
    Present    &  -1.7482,  0.4253,   -6.7611 &  8\\
    Ghosh and Roy \cite{ghosh2009frame} & -1.7554,  0.4156,  -6.7485  & 10\\
    Ibrahimbegovic and Taylor \cite{ibrahimbegovic2002role} &  -1.7368,  0.4308,  -6.7329  & 8 \\
    Jeleni\'c and Crisfield \cite{jelenic1999geometrically} & \  \  \  \   \  \ - \  \ , \  \  \  \   - \  \ , -6.7684 & 2\\
    \hline
   \end{tabular}
   \label{tab:frameinvariance_table}
   \label{tab:LCantilever}
\end{table}

\subsection{Performance study: Cantilever subjected to torsion and moment}
\begin{figure} [h]
  \begin{subfigure}{.45\textwidth}
    \centering
    \includegraphics[scale=0.40]{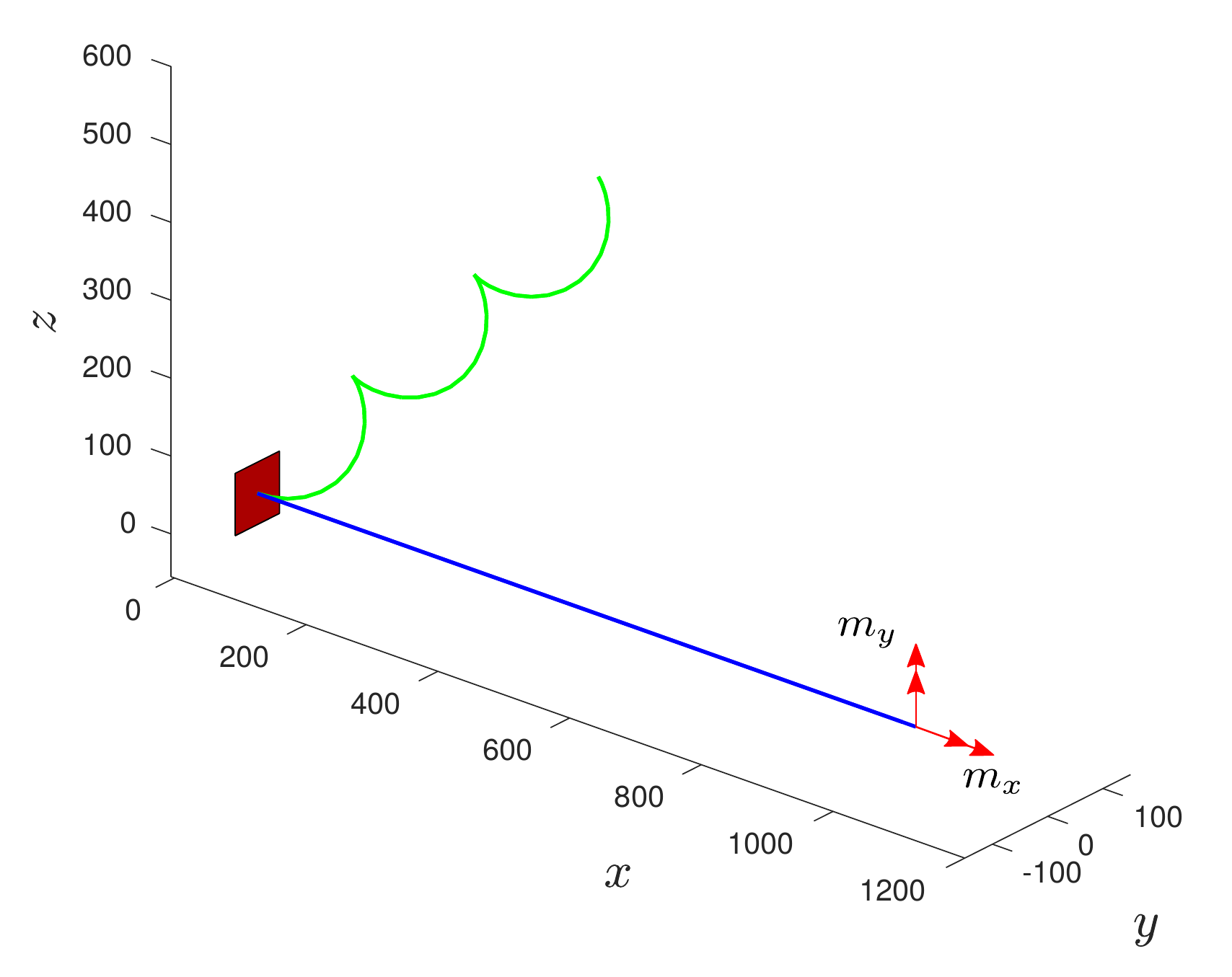}
    \caption{}
 \label{fig:spiralBendingTorsion1}
  \end{subfigure}
  \begin{subfigure}{.45\textwidth}
    \centering
    \includegraphics[scale=0.31]{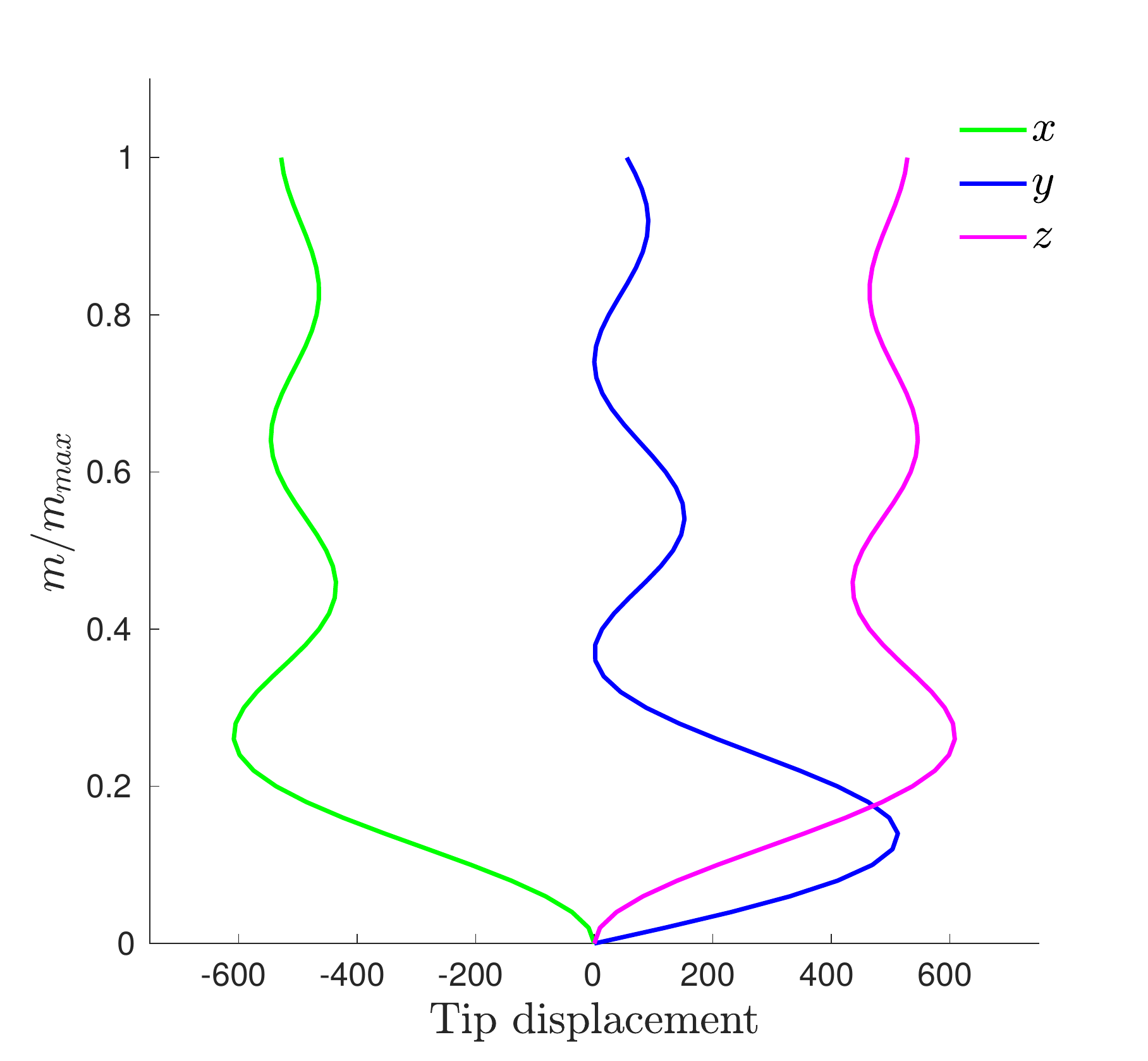}
    \caption{}
 \label{fig:spiralBendingTorsion2}
  \end{subfigure}
  \caption{(a) Reference and deformed configurations of the cantilever beam  when subjected to applied end moments; (b) load versus tip displacement predicted by our present mixed FE model}
  \label{fig:spiralBendingTorsion}
  \end{figure}

The performance of the mixed FE model is evaluated by the computational time required to obtain a solution,  which in turn depends on the type of formulation, size of the stiffness matrix and number of Newton iterations required for convergence. In this subsection, we study the number of Newton iterations and compare it with other
nonlinear beam finite elements found in the literature.  To demonstrate the
performance of the present formulation, we study a straight  cantilever beam
subjected to end moment and torque. This problem was considered by
\cite{schulz2019finite} in the context of number of Newton iterations, while
\cite{meier2015locking} studied it for its extremely large rotation. The
cantilever beam has a length 1000 units. The straight reference
configuration and the boundary conditions are shown in  Figure
\ref{fig:spiralBendingTorsion1}. The stiffness constants are $EI_1= EI_2 = GJ =8.333\times10^2$  and $AE=100$. The initial configuration of the beam is discretized using 100 linear finite elements. Under the action of external moments, the beam undergoes complex 3D rotations and takes the shape of a spiral aligned at
$45^\circ$ from the $x-$axis. Figure \ref{fig:spiralBendingTorsion2} shows the
deformed configurations predicted by the proposed numerical procedure. The
maximum moment of $m_x=10$ and $m_y=10$ is reached in a single load step. The
number of Newton iterations required to obtain results is compared with other
methods as shown in Table  \ref{tab:spiralPerformance}. From the table, it can
be seen that the present formulation converges very fast  even in the case of extremely large deformation. Moreover, Table \ref{tab:spiralPerformance} also shows that the number of Newton iteration required for the convergence of the solution is independent of the finite element mesh. 

\begin{table}
  \centering
  \caption{Comparison of number of Newton iterations with number of elements used for discretization}
  \begin{tabular}{l| c| c}
    \hline
    Authors &  \vtop{\hbox{\strut \# Newton iteration}\hbox{\strut \;\;\;\;\;\;\;(8 elements)}} & \vtop{\hbox{\strut \# Newton iteration}\hbox{\strut \;\;\;\;\;\;\;(16 elements)}}\\
    \hline
    Present    & 3    & 3 \\
    Schulz and B\"ol \cite{schulz2019finite}&  15  & 15 \\
   Meier \textit{et al.} \cite{meier2014objective} &  438  & 625 \\
   Meier \textit{et al.} \cite{meier2019geometrically} &875 & 1220\\
   Simo \cite{simo1985finite} & 892  & 1011 \\
   Jeleni\'c and Crisfield \cite{jelenic1999geometrically} &1497 & 1540\\
    \hline
   \end{tabular}
   \label{tab:spiralPerformance}
\end{table}

\subsection{Star shaped dome}
\begin{figure}
    \centering
    \includegraphics[scale=0.5]{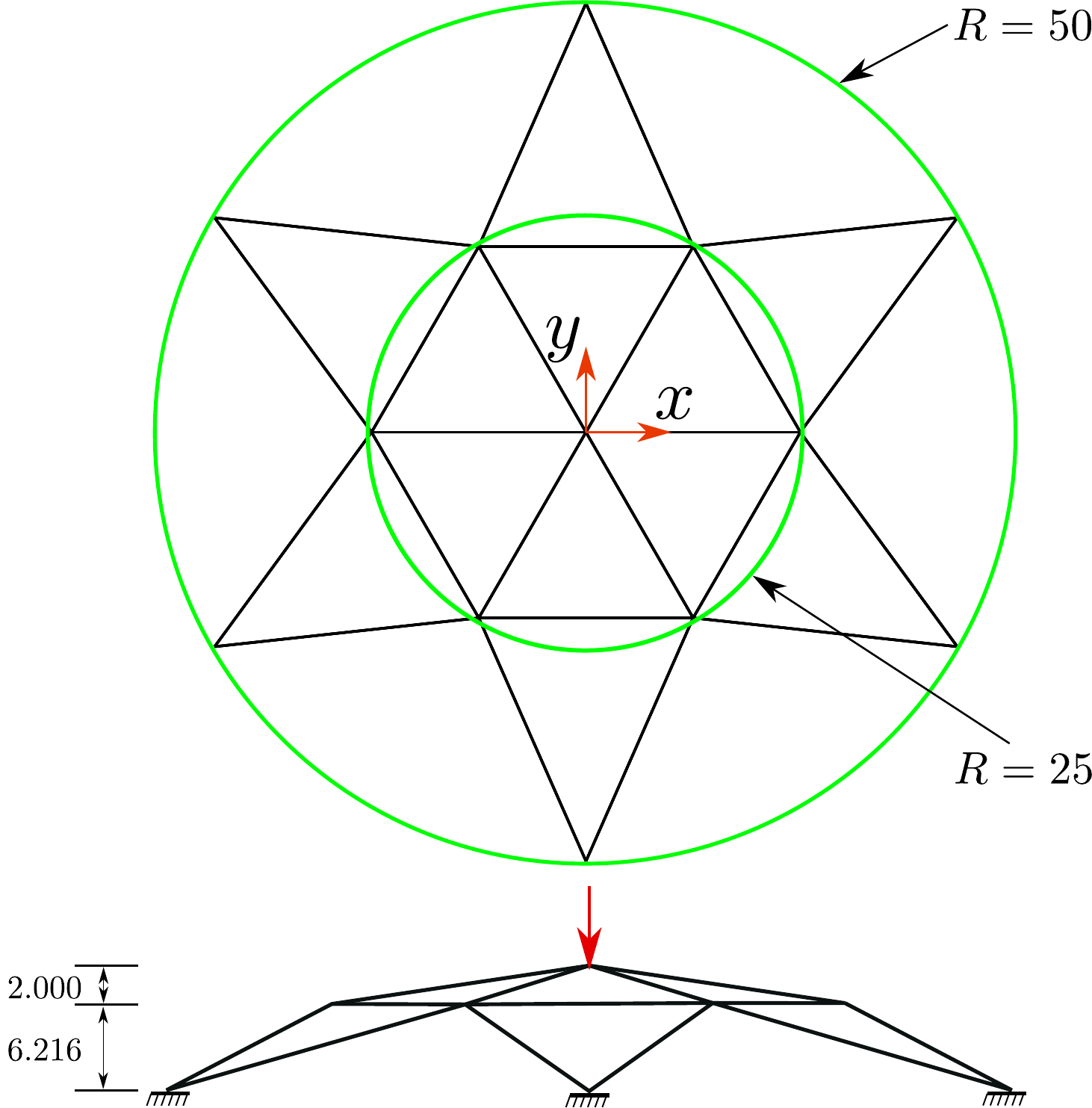}
    \caption{Geometric details of space dome with point load applied at the apex}
  \label{fig:starDome0}
  \end{figure}

\begin{figure}
    \begin{subfigure}{.45\textwidth}
      \centering
      \includegraphics[scale=0.85]{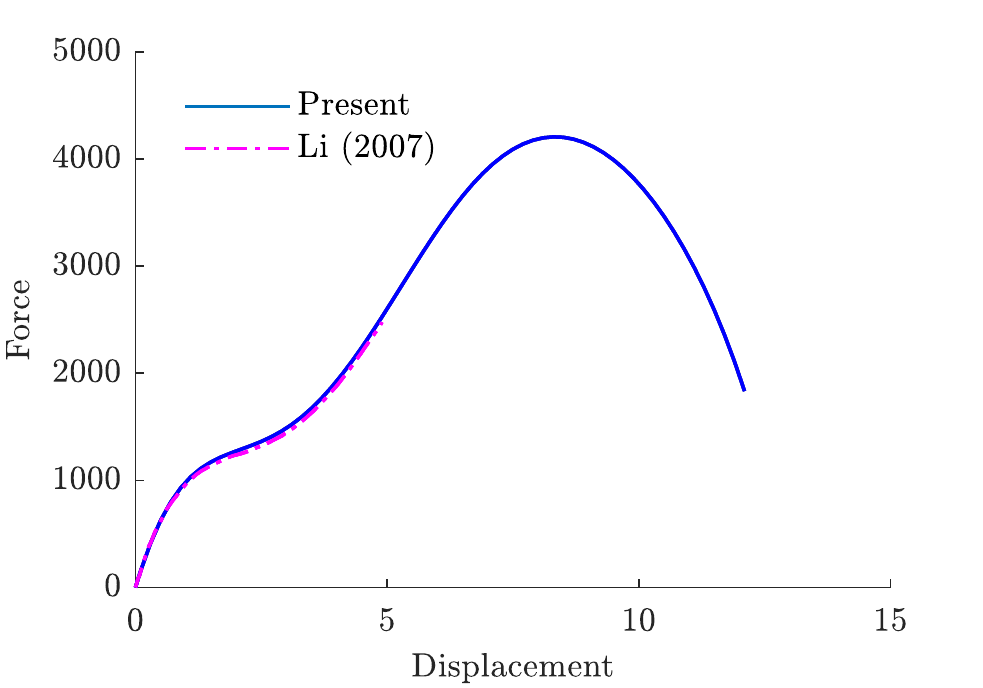}
      \caption{}
    \label{fig:starDome1}
    \end{subfigure}%
    \begin{subfigure}{.45\textwidth}
      \centering
      \includegraphics[scale=0.85]{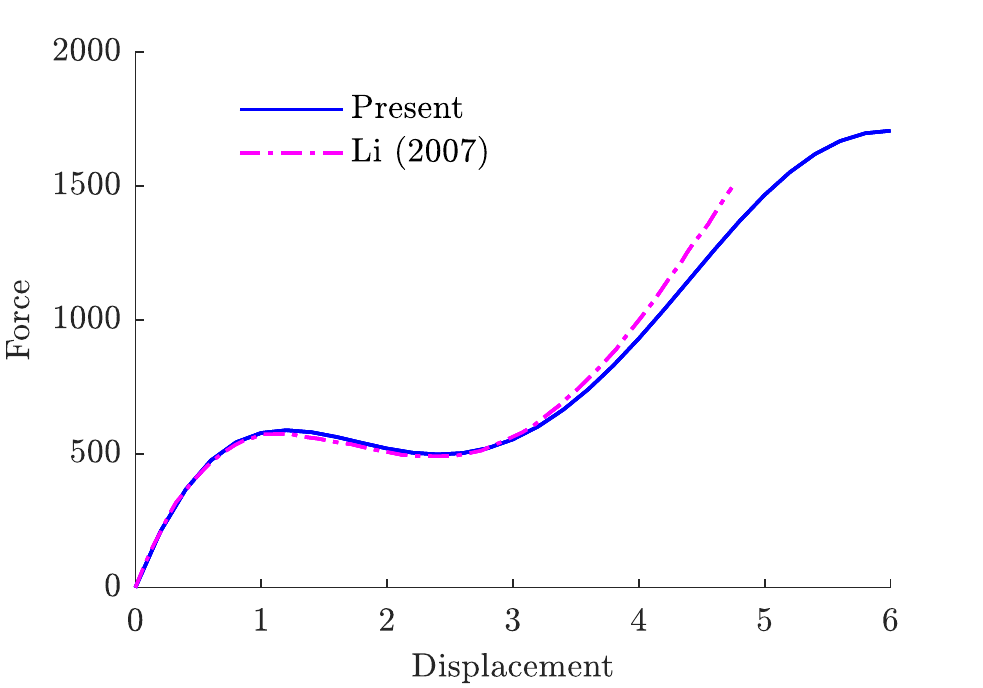}
      \caption{}
    \label{fig:starDome2}
    \end{subfigure}
    \caption{(a) Load-displacement plot for square cross-section; (b)  load-displacement plot for rectangular cross-section }
    \label{fig:starDome}
  \end{figure}

Finally, we study the post-buckling behavior of a 3D frame structure. This problem is chosen to demonstrate the
ability of our present numerical technique to handle multiple member
intersection at a common node. Displacement based FE techniques for Kirchhoff
 rods require higher order smoothness across the element boundary. For
example, techniques based on isogeometric analysis have higher continuity across a patch of finite elements; additional constraints are invoked to impose the required continuity at the intersection of multiple members. The present FE
technique completely avoids these additional constraints. Moreover, multiple
member intersection may be effectively handled within the finite element assembly. We demonstrate these properties of our FE formulation by applying it to study the deformation of a 24 member star shaped shallow dome. The geometry and boundary conditions of the star shaped dome are shown in figure \ref{fig:starDome0}. 
A point load is applied at the crown and the simulation is performed using  two sets of cross sectional  properties.  For the first case, we take $AE=9.6051\times10^5$, $GJ=1.5465\times10^5$, and  $EI_1=EI_2=2.5361\times10^5$. The load displacement curve predicted by our FE technique is presented in Figure \ref{fig:starDome1}. For a comparison, the load displacement curve predicted by Li \cite{li2007co} for the same stiffness properties is also shown in Figure \ref{fig:starDome1}. Similarly, the load displacement curve predicted by our FE technique for material stiffness constants $ GJ= 1.0061\times10^5$, $EI_1= 0.8938 \times 10^5$, and $EI_2=7.2023\times10^5 $ is shown in Figure \ref{fig:starDome2}. For a comparison, the load displacement computed by Li \cite{li2007co} is also reproduced. These two simulations make it amply clear that our FE approximation is capable of predicting post-buckling behaviour even in the case of multiple members.


\section{Summary}
\label{sec:conclusion}
Cast in the language of exterior calculus, we formulated a four-field mixed  variational
principle for Kirchhoff rods in this article. The variational principle was general enough to accommodate major modes of deformation (bending, twisting and extension) in a three dimensional setting. Conceptually, the mixed variational principle can be understood as a rule to evolve a frame bundle,  whose base space is the rod's centerline (a 1-dimensional manifold). This principle utilised Cartan's moving frame to describe the geometric properties of the rod configurations.  Building on this theoretical construction, we presented a mixed FE model for Kirchhoff rods, that is only $C^0$. From the perspective of the moving frame, the present mixed FE approximation may be interpreted as a discrete approximation of the frame bundle. The
interpolation schemes used to approximate the frame fields were chosen
to respect the underlying geometric structure. These features of our
approximation scheme lead to a remarkably superior performance of the numerical method with
only a marginal increase in the number of DOF's. 

\section*{Acknowledgements}
JK was supported by ISRO through the Centre of Excellence in Advanced Mechanics of Materials grant No. ISRO/DR/0133. BD and JNR gratefully acknowledges the support by DoD-Army Research grant (DARES) W911NF2120064.
\newpage
\bibliographystyle{abbrv}
\bibliography{kirchhoff}
\appendix
\section{Rotation group}
\label{section2:Rotation group}
The rotation group plays a central role in the variational principle and the
numerical procedure we developed in the previous sections. For completeness, we present a brief introduction to the rotation group. The set of all
orientation preserving rotations in $\mathbb{R}^3$  is denoted by $SO(3)$. For
a vector $\bs{v}\in \mathbb{R}^3$, the action of $SO(3)$ is given by,
\[
  \bs{v}'=\bs{\Lambda} \bs{v}; \quad \bs{\Lambda}\in SO(3),
\]
where $\bs{v}'$ is the rotated vector. In the above equation, the action of $\Lambda$
on $\bs{v}$ is linear. An important property of $SO(3)$ is that it preserves the
inner-product, for two vectors $\bs{v},\bs{w} \in \mathbb{R}^3$, we have $\langle \bs{v},\bs{w}
\rangle_{\mathbb{R}^3}=\langle \bs{\Lambda} \bs{v},\bs{w} \Lambda \rangle_{\mathbb{R}^3}$,
where $\langle.,. \rangle_{\mathbb{R}^3}$ is the inner-product on
$\mathbb{R}^3$. Moreover, it is easy to see that $\norm{\bs{\Lambda
v}}_{\mathbb{R}^3}=\norm{\bs{v}}_{\mathbb{R}^3}$, where  $\norm{.}_{\mathbb{R}^3}$
denotes the Euclidean norm.  Using  the above properties, the set of all
rotations can be formally defined as,
\begin{equation}
  SO(3)=\{\Lambda|\;
  \norm{\Lambda \bs{v}}_{\mathbb{R}^3}=\norm{\bs{v}}_{\mathbb{R}^3}\; \forall \bs{v}\in \mathbb{R}^3; \det(\bs{\Lambda})=1 \}.
\end{equation}
From the definition of $SO(3)$, it is easy to see that for two rotations $\bs{\Lambda}_1$ and $\bs{\Lambda}_2$, their composition
$\bs{\Lambda}_1\circ\bs{\Lambda}_2$ is also a rotation. More generally, the set $SO(3)$ is
closed under composition, making it a group with composition as the group operation. Since each rotation is a linear map  on
$\mathbb{R}^3$, we can identify elements form  $SO(3)$ with $3 \times 3$
matrices. Moreover, the restriction that elements from $SO(3)$ preserves the
length of vectors permits us to identify the rotation group with a set of
orthogonal matrices (a subset of $3 \times 3$  matrices). 
Using the
matrix representation, composition of two rotations can be written as a matrix
product of the associated rotation matrices. In the previous discussion, we did
not distinguish between the rotation group and its matrix representation.
Indeed, we use $SO(3)$ to denote the set of all orthogonal matrices with
determinant one; we do the same in the following sections as well.

From the matrix representation of $SO(3)$, it looks like nine components are
required to specify a rotation. However the inter-relationships dictated by
$\bs{\Lambda}^t\bs{\Lambda}=\bs{I}$ and $\det \bs{\Lambda}=1$, reduces the number of independent
components to just three. Another important property of $SO(3)$ is, it is
also a smooth manifold. On $SO(3)$, one can smoothly transform one rotation to
another; in other words, one can connect two elements in $SO(3)$ using a smooth
curve. These smooth curves permit us to talk about the tangent vector to $SO(3)$. Combining the facts that $SO(3)$ is a smooth manifold and the number of
independent parameters required to specify a rotation is three, we arrive at the
dimension of $SO(3)$ as three. The dimension of a smooth manifold tells us the
number of independent coordinates required to describe a small smooth patch on
the manifold. For $SO(3)$, Euler angles, rotation vector are examples of
parametrization with three parameters. It should be mentioned that none of these
parametrization is global; in other words, two or more coordinate covers are
required to parametrize the entire of $SO(3)$.

\subsection{Tangent space of $SO(3)$}
The tangent space to a manifold is a local linear approximation to the manifold.
Important kinematic quantities like velocity can be identified as an element
from the tangent space of the configuration manifold. In this subsection, we
discuss the tangent space of $SO(3)$. A geometric approach to the tangent space of a manifold is to look at the tangent vectors of curves passing through a point. Consider a smooth curve $\gamma:(0,1)\rightarrow
SO(3)$. Along the curve $\gamma$, we have the following relationship,
\[
  \bs{\Lambda}^t(s)\bs{\Lambda}(s)=\bs{I};\quad s\in(0,1)
  \]
Now differentiating the above equation with respect to $s$, we have,
\begin{equation}
  \frac{d\bs{\Lambda}^t}{ds}\bs{\Lambda}+\bs{\Lambda}^t\frac{d\bs{\Lambda}}{ds}=0.
  \label{eq:derOrthonormal}
\end{equation}
Using the properties of matrix multiplication and transpose, we have $\frac{d\bs{\Lambda}^t}{ds}\bs{\Lambda}=\left[\bs{\Lambda}^t\frac{d\bs{\Lambda}}{ds}\right]^t$, using this relationship in the above equation, we conclude that $\frac{d\bs{\Lambda}^t}{ds}\bs{\Lambda}$ is skew symmetric. We denote the set of all skew symmetric matrices with dimension three as $so(3)$. 
Using the hat map, $so(3)$ can be identified with $\mathbb{R}^3$.  We define $\widehat{\bs{\Omega}}:=-\frac{d\bs{\Lambda}^t}{ds}\bs{\Lambda}$, where $\bs{\Omega}\in\mathbb{R}^3$ and $\widehat{(.)}$ is the hat map which sends $\mathbb{R}^3$ to $so(3)$. By using the definition of $\widehat{\bs{\Omega}}$ in  \eqref{eq:derOrthonormal}, we arrive at
\[
    \frac{d \Lambda}{ds}=\Lambda \widehat{\bs{\Omega}}.
\]
Using the definition of $\gamma$, the above equation can be rewritten as,
\begin{equation}
  \frac{d \gamma}{ds}=\bs{\Lambda} \widehat{\bs{\Omega}};\quad \gamma(s)=\bs{\Lambda}\in SO(3),\; \widehat{\bs{\Omega}}\in so(3).
\end{equation}
The above equation gives an expression for the tangent vector to the curve $\gamma$. By  changing the curve $\gamma$ suitably, one can vary the $so(3)$ factor in the expression for tangent vector. Thus the tangent space of $SO(3)$ at $\bs{\Lambda}\in SO(3)$ can be formally defined as,
\begin{equation}
    T_{\bs{\Lambda}}SO(3)=\{\bs{\Lambda}\widehat{\bs{\Omega}}| \widehat{\omega}\in so(3)\}.
    \label{eq:rightTranslate}
\end{equation}
Similarly, one could start from the expression $\bs{\Lambda\Lambda}^t=\bs{I}$ and repeat the same steps and arrive at the following expression for  $T_{\bs{\Lambda}}SO(3)$
\begin{equation}
  T_{\bs{\Lambda}}SO(3)=\{\widehat{\bs{\omega}}\bs{\Lambda}| \widehat{\bs{\omega}}\in so(3)\}.
  \label{eq:leftTranslate}
\end{equation}
The linear space defined in  \eqref{eq:rightTranslate} and \eqref{eq:leftTranslate} are exactly same; these tangent spaces are called the left and right translations of $SO(3)$. The $so(3)$ components of the right and left translated tangent spaces are related by the following expression,
\begin{equation}
  \widehat{\bs{\Omega}}=\bs{\Lambda} \widehat{\bs{\omega}}\bs{\Lambda}^t.
\end{equation}
We now produce a basis for $so(3)$. The canonical basis for $\mathbb{R}^3$ is denoted by $\{\bs{w}_1,\ldots,\bs{w}_3\}$. Since the hat-map is an invertible linear map, the action of the hat-map on the canonical basis for $\mathbb{R}^3$ produces a basis for $so(3)$. This basis for $so(3)$ is denoted as $\mathcal{B}_{so(3)}=\{\widehat{\bs{w}}_1, \widehat{\bs{w}}_2, \widehat{\bs{w}}_3\}$. The basis elements of $so(3)$ obtained by applying the hat-map to the canonical basis of $\mathbb{R}^3$ is given as,
\begin{eqnarray}
    \widehat{\bs{w}}_1=\begin{bmatrix}   0 & 0 & 0\\  0 & 0 & -1\\ 0 & 1 & 0 \end{bmatrix};\quad
    \widehat{\bs{w}}_2=\begin{bmatrix}    0 & 0 & 1\\ 0 & 0 & 0\\-1 & 0 & 0 \end{bmatrix};\quad
    \widehat{\bs{w}}_3= \begin{bmatrix}  0 & -1 & 0 \\ 1 & 0 & 0\\ 0 & 0 & 0 \end{bmatrix}.
\end{eqnarray}
In terms of the basis $\mathcal{B}_{so(3)}$, any element $\widehat{\bs{w}}\in so(3)$ may now  be written as
\begin{eqnarray}
  \widehat{\bs{w}}=w_1\widehat{\bs{w}}_1+w_2\widehat{\bs{w}}_2+w_3\widehat{\bs{w}}_3,
  \end{eqnarray}
where $w_i$ are real numbers. Thus the skew symmetric matrix associated with $\widehat{\bs{w}}$ can be written as,
\begin{eqnarray}
  \widehat{\bs{w}} = \begin{bmatrix}
  0 & -w_3 & w_2\\
  w_3 & 0 & -w_1\\
  -w_2 & w_1 & 0
\end{bmatrix}.
\end{eqnarray}
The axial vector of an element from $so(3)$ is obtained by applying the inverse of the hat-map. The inverse hat-map denoted by $\widetilde{(\cdot)}$, sends an element from $so(3)$ to $\mathbb{R}^3$. Now applying the inverse hat-map to the $\widehat{w}$, we arrive at $\widetilde{\widehat{\bs{w}}}=w_1\bs{w}_1+w_2\bs{w}_2+w_3\bs{w}_3$.

The Lie bracket on $T_{\Lambda}SO(3)$ is a bilinear map which takes $T_{\Lambda}SO(3)$ to itself, it is denoted by $[\cdot,\cdot]$. The action of Lie bracket on the basis vectors of $T_ISO(3)$ can be written as,
\begin{equation}
  [\widehat{\bs{w}}_1, \widehat{\bs{w}}_2] = \widehat{\bs{w}}_3, \quad
  [\widehat{\bs{w}}_3, \widehat{\bs{w}}_1] = \widehat{\bs{w}}_2, \quad
  [\widehat{\bs{w}}_2, \widehat{\bs{w}}_3] = \widehat{\bs{w}}_1
\end{equation}
Using bilinearity, the above definition can be extended to any to elements from $T_I SO(3)$. In terms of the axial vector, the Lie bracket can be written as,
\begin{equation}
  \bs{w}_1\times \bs{w}_2=\bs{w}_3;\quad \bs{w}_2\times \bs{w}_3=\bs{w}_1;\quad \bs{w}_3\times \bs{w}_1=\bs{w}_2.
\end{equation}
where $\times$ denotes the cross product between two vectors in $\mathbb{R}^3$.

Having introduced the tangent space on $SO(3)$, we now introduce the notion of infinitesimal length or metric on $SO(3)$. Metric is a positive-definite, symmetric, bilinear map from the tangent space to $\mathbb{R}$. On $SO(3)$, we introduce the following metric:
\begin{equation}
  \left\langle \bs{u},\bs{v} \right\rangle_{T_{\Lambda} SO(3)}:=\frac{1}{2}\tr(\bs{u}^t\bs{v}),
  \label{eq:innerproductTSO(3)}
\end{equation}
for $\bs{u},\bs{v} \in T_\Lambda SO(3)$. The above metric is bi-invariant, that is the metric produces the same length for both right and left translations of $so(3)$. The metric on $T_I SO(3)$ can be pulled back to $\mathbb{R}^3$ using the inverse hat-map. This relationship can be written as,
\begin{equation}
  \widetilde{}\;{}^*(\left\langle \bs{u},\bs{v} \right\rangle_{T_{I} SO(3)})=\langle \widetilde{\bs{u}},\widetilde{\bs{v}} \rangle_{\mathbb{R}^3},
\end{equation}
where $\widetilde{}\;{}^*(.)$ denotes the pull-back under the inverse hat-map. While the inner-product  on the left is defined on $T_{I} SO(3)$, the right-hand side is an Euclidean inner-product between vectors in $\mathbb{R}^3$. Since the metric is bi-invariant, the pull-back  map can be extended to any tangent space of $SO(3)$. 

\subsection{Affine connection on $SO(3)$}
The bi-invariant metric introduced on $SO(3)$ makes it a Riemannian manifold.
This metric produces a unique connection which can be used to parallel
transport tangent vectors. The connection on $SO(3)$ also provides a coordinate
independent differential operator called the covariant derivative. For vector fields $u$ and $v$, the covariant derivative $v$ in the direction $u$ is given by
\begin{equation}
  \nabla_{\widehat{\bs{u}}}\widehat{\bs{v}}=\frac{1}{2}[\widehat{\bs{u}},\widehat{\bs{v}}].
\end{equation}
In terms of the  associated axial vectors the above relationship may be written
as,
\begin{equation}
  \langle \nabla_{\widehat{\bs{u}}}  \widehat{\bs{v}},\widehat{\bs{w}}  \rangle_{T_{\Lambda}SO(3)}=\frac{1}{2}\langle (\bs{u} \times \bs{v}), \bs{w}\rangle_{\mathbb{R}^3}.
\end{equation}
The covariant derivative plays a key role in the computation of Hessian in the numerical solution procedure; for more details about the role of connection, the reader may refer to Simo \cite{simo1992symmetric}.


\section{Residual and tangent calculation}
\label{sec:residueTangent}
For the purpose of computing the residual and tangent vectors, we ignore the differences between 1-forms and vectors and treat both of them as vectors; this leads to a description of stiffness coefficients very similar to the one encountered in conventional displacement-based finite element model. The components of the residual vector can be computed as, 
\begin{align}
    \nonumber
    \ms{D}_{\bs{\eta}} I^h&=\int_{[0,L^h]} AE(\eta-1)-n^1\\
    \nonumber
    \ms{D}_{\bs{n}^1} I^h&=\int_{[0,L^h]} -\eta+ \ed \bs{u}^t ( \bs{\Lambda\Lambda}_0) \bar{\tb{E}}_1 ;\quad
    \ms{D}_{\bs{n}^2} I^h=\int_{[0,L^h]}  \ed \bs{u}^t ( \bs{\Lambda\Lambda}_0) \bar{\tb{E}}_2 ;\quad
    \ms{D}_{\bs{n}^3} I^h=\int_{[0,L^h]}  \ed \bs{u}^t (\bs{\Lambda\Lambda}_0) \bar{\tb{E}}_3 \\
    \nonumber
    \ms{D}_{\tb u^1} I^h&=\int_{[0,L^h]} ({\bar{\tb{E}}_1}^t ( \bs{\Lambda\Lambda}_0) \bs{n} ) \ed \tb N ;\quad
    \ms{D}_{\tb u^2} I^h=\int_{[0,L^h]} ({\bar{\tb{E}}_2}^t ( \bs{\Lambda\Lambda}_0) \bs{n} ) \ed \tb N ;\quad
    \ms{D}_{\tb u^3} I^h=\int_{[0,L^h]} ({\bar{\tb{E}}_3}^t ( \bs{\Lambda\Lambda}_0) \bs{n} ) \ed \tb N \\
    \nonumber
  \ms{D}_{\bs \lambda^1} I^h&=\int_{[0,L^h]}  \ms D_{\bs \lambda^1} W_b^h+ \ed \bs{\bs{u}}^t \ms{D}_{\bs \lambda^1} ( \bs{\Lambda\Lambda}_0) \bs{n};\quad
  \ms{D}_{\bs \lambda^2} I^h=\int_{[0,L^h]}  \ms D_{\bs \lambda^2} W_b^h+ \ed \bs{u}^t \ms{D}_{\bs \lambda^2} ( \bs{\Lambda\Lambda}_0) \bs{n}\\
    \nonumber
  \ms{D}_{\bs \lambda^3} I^h&=\int_{[0,L^h]}  \ms D_{\bs \lambda^3} W_b^h+ \ed \bs{u}^t \ms{D}_{\bs \lambda^3} ( \bs{\Lambda\Lambda}_0) \bs{n}
\end{align}
The components of tangent operator can be expressed as,
\begin{align}
 \nonumber
 &\ms D_{\bs \eta}\ms D_{\bs \eta}I^h=\int_{[0,L]}  AE;\quad
 \ms D_{\bs \eta}\ms D_{\tb t^1}I^h=\int_{[0,L]} -1 \\
\nonumber
 &\ms D_{\tb n^1 \bs \lambda^1}I^h=\int_{[0,L^h]}  (\ed \bs{u}^t \ms{D}_{\bs \lambda^1} ( \bs{\Lambda\Lambda}_0) \bar{\tb{E}}_1 );\quad
 \ms D_{\tb n^1\bs \lambda^2}I^h=\int_{[0,L^h]} ( \ed \bs{u}^t \ms{D}_{\bs \lambda^2} ( \bs{\Lambda\Lambda}_0) \bar{\tb{E}}_1 );\quad
 \ms D_{\tb n^1\bs \lambda^3}I^h=\int_{[0,L^h]}  (\ed \bs{u}^t \ms{D}_{\bs \lambda^3} ( \bs{\Lambda\Lambda}_0) \bar{\tb{E}}_1 ) \\ 
\nonumber
 &\ms D_{\tb n^2 \bs \lambda^1}I^h=\int_{[0,L^h]}  (\ed \bs{u}^t \ms{D}_{\bs \lambda^1} ( \bs{\Lambda\Lambda}_0) \bar{\tb{E}}_2 );\quad
 \ms D_{\tb n^2\bs \lambda^2}I^h=\int_{[0,L^h]} ( \ed \bs{u}^t \ms{D}_{\bs \lambda^2} ( \bs{\Lambda\Lambda}_0) \bar{\tb{E}}_2 )\quad
 \ms D_{\tb n^2\bs \lambda^3}I^h=\int_{[0,L^h]}  (\ed \bs{u}^t \ms{D}_{\bs \lambda^3} ( \bs{\Lambda\Lambda}_0) \bar{\tb{E}}_2 ) \\ 
\nonumber
 &\ms D_{\tb n^3\bs \lambda^1}I^h=\int_{[0,L^h]}  (\ed \bs{u}^t \ms{D}_{\bs \lambda^1} ( \bs{\Lambda\Lambda}_0) \bar{\tb{E}}_3 );\quad
 \ms D_{\tb n^3\bs \lambda^2}I^h=\int_{[0,L^h]} ( \ed \bs{u}^t \ms{D}_{\bs \lambda^2} ( \bs{\Lambda\Lambda}_0) \bar{\tb{E}}_3 );\quad
 \ms D_{\tb n^3\bs \lambda^3}I^h=\int_{[0,L^h]}  (\ed \bs{u}^t \ms{D}_{\bs \lambda^3} ( \bs{\Lambda\Lambda}_0) \bar{\tb{E}}_3 )\\
\nonumber
 &\ms D_{\tb n^1\tb u^1}I^h=\int_{[0,L^h]}  ({\bar{\tb{E}}_1}^t ( \bs{\Lambda\Lambda}_0) \bar{\tb{E}}_1 )\ed \tb N;\quad
 \ms D_{\tb n^1\tb u^2}I^h=\int_{[0,L^h]}  ({\bar{\tb{E}}_2}^t ( \bs{\Lambda\Lambda}_0) \bar{\tb{E}}_1 )\ed \tb N;\quad
 \ms D_{\tb n^1\tb u^3}I^h=\int_{[0,L^h]}  ({\bar{\tb{E}}_3}^t ( \bs{\Lambda\Lambda}_0) \bar{\tb{E}}_1 )\ed \tb N \\ 
\nonumber
 &\ms D_{\tb n^2\tb u^1}I^h=\int_{[0,L^h]}  ({\bar{\tb{E}}_1}^t ( \bs{\Lambda\Lambda}_0) \bar{\tb{E}}_2 )\ed \tb N;\quad
 \ms D_{\tb n^2\tb u^2}I^h=\int_{[0,L^h]}  ({\bar{\tb{E}}_2}^t ( \bs{\Lambda\Lambda}_0) \bar{\tb{E}}_2 )\ed \tb N;\quad
 \ms D_{\tb n^2\tb u^3}I^h=\int_{[0,L^h]}  ({\bar{\tb{E}}_3}^t ( \bs{\Lambda\Lambda}_0) \bar{\tb{E}}_2 )\ed \tb N \\ 
\nonumber
 &\ms D_{\tb n^3\tb u^1}I^h=\int_{[0,L^h]}  ({\bar{\tb{E}}_1}^t ( \bs{\Lambda\Lambda}_0) \bar{\tb{E}}_3 )\ed \tb N;\quad
 \ms D_{\tb n^3\tb u^2}I^h=\int_{[0,L^h]}  ({\bar{\tb{E}}_2}^t ( \bs{\Lambda\Lambda}_0) \bar{\tb{E}}_3 )\ed \tb N;\quad
 \ms D_{\tb n^3\tb u^3}I^h=\int_{[0,L^h]}  ({\bar{\tb{E}}_3}^t ( \bs{\Lambda\Lambda}_0) \bar{\tb{E}}_3 )\ed \tb N\\
\nonumber
 &\ms D_{\tb u^1\bs \lambda^1}I^h=\int_{[0,L^h]}  \ed \tb N   \otimes (\bar{\tb{E}}_1^t \ms{D}_{\bs \lambda^1} ( \bs{\Lambda\Lambda}_0) \bs{n});\quad
 \ms D_{\tb u^1\bs \lambda^2}I^h=\int_{[0,L^h]} \ed \tb N   \otimes (\bar{\tb{E}}_1^t \ms{D}_{\bs \lambda^2} ( \bs{\Lambda\Lambda}_0) \bs{n}) \\
\nonumber
&\ms D_{\tb u^1\bs \lambda^3}I^h=\int_{[0,L^h]}  \ed \tb N   \otimes(\bar{\tb{E}}_1^t \ms{D}_{\bs \lambda^3} ( \bs{\Lambda\Lambda}_0) \bs{n}) ;\quad
 \ms D_{\tb u^2\bs \lambda^1}I^h=\int_{[0,L^h]} \ed \tb N   \otimes  (\bar{\tb{E}}_2^t \ms{D}_{\bs \lambda^1} ( \bs{\Lambda\Lambda}_0) \bs{n})\\
\nonumber
 &\ms D_{\tb u^2\bs \lambda^2}I^h=\int_{[0,L^h]} \ed \tb N   \otimes (\bar{\tb{E}}_2^t \ms{D}_{\bs \lambda^2} ( \Lambda\Lambda_0) \bs{n});\quad
 \ms D_{\tb u^2\bs \lambda^3}I^h=\int_{[0,L^h]} \ed \tb N   \otimes (\bar{\tb{E}}_2^t \ms{D}_{\bs \lambda^3} ( \bs{\Lambda\Lambda}_0)\bs{n}) \\
\nonumber
 &\ms D_{\tb u^3\bs \lambda^1}I^h=\int_{[0,L^h]} \ed \tb N   \otimes  (\bar{\tb{E}}_3^t \ms{D}_{\bs \lambda^1} ( \bs{\Lambda\Lambda}_0) \bs{n});\quad
 \ms D_{\tb u^3\bs \lambda^2}I^h=\int_{[0,L^h]}  \ed \tb N   \otimes (\bar{\tb{E}}_3^t \ms{D}_{\bs \lambda^2} ( \bs{\Lambda\Lambda}_0) \bs{n})\\
\nonumber
 &\ms D_{\tb u^3\bs \lambda^3}I^h=\int_{[0,L]}  \ed \tb N   \otimes (\bar{\tb{E}}_3^t \ms{D}_{\bs \lambda^3} ( \bs{\Lambda\Lambda}_0) \bs{n})\\
\nonumber
 &\ms D_{\bs \lambda^1\bs \lambda^1}I^h=\int_{[0,L]}  \ms D_{\bs \lambda^1\bs \lambda^1} W_b^h + (\ed \bs{u}^t  \ms{D}_{\bs \lambda^1\bs \lambda^1} ( \bs{\Lambda\Lambda}_0) \bs{n}); \quad \ms D_{\bs \lambda^2\bs \lambda^2}I^h=\int_{[0,L]}  \ms D_{\bs \lambda^2\bs \lambda^2} W_b^h+  (\ed u^t  \ms{D}_{\bs \lambda^2\bs \lambda^2} ( \bs{\Lambda\Lambda}_0) \bs{n})\\
\nonumber
 &\ms D_{\bs \lambda^3\bs \lambda^3}I^h=\int_{[0,L]}  \ms D_{\bs \lambda^3\bs \lambda^3} W_b^h+  (\ed \bs{u}^t  \ms{D}_{\bs \lambda^3\bs \lambda^3} ( \bs{\Lambda\Lambda}_0) \bs{n}); \quad \ms D_{\bs \lambda^1\bs \lambda^2}I^h=\int_{[0,L]}  \ms D_{\bs \lambda^1\bs \lambda^2} W_b^h+  (\ed \bs{u}^t \ms{D}_{\bs \lambda^1\bs \lambda^2} ( \bs{\Lambda\Lambda}_0) \bs{n})\\
\nonumber
 &\ms D_{\bs \lambda^1\bs \lambda^3}I^h=\int_{[0,L]}  \ms D_{\bs \lambda^1\bs \lambda^3} W_b^h+  (\ed \bs{u}^t  \ms{D}_{\bs \lambda^1\bs \lambda^3} ( \bs{\Lambda\Lambda}_0) \bs{n}); \quad 
 \ms D_{\bs \lambda^2\bs \lambda^3}I^h=\int_{[0,L]}  \ms D_{\bs \lambda^2\bs \lambda^3} W_b^h+  (\ed \bs{u}^t \ms{D}_{\bs \lambda^2\bs \lambda^3} ( \bs{\Lambda\Lambda}_0) \bs{n})
\end{align}
Where 
\begin{align}
\nonumber
\ms D_{\bs \lambda^i}(\bs{\Lambda\Lambda}_0)&=  \widehat{\tb{E}}_i(\bs{\Lambda\Lambda}_0) \otimes \tb N, \\
\nonumber
\ms D_{\bs \lambda^i\bs \lambda^j}(\bs{\Lambda\Lambda}_0)&=\frac{1}{2}( \widehat{\tb{E}}_i \widehat{\tb{E}}_j+ \widehat{\tb{E}}_j \widehat{\tb{E}}_i) (\bs{\Lambda\Lambda}_0)\otimes(\tb N \otimes \tb N).
\end{align}
\end{document}